\documentclass[12pt]{article}

\usepackage{amsmath,amssymb,amsthm}
\usepackage[figuresright]{rotating} 
\usepackage{mathrsfs,array,booktabs}
\usepackage[authoryear]{natbib}
\usepackage{hyperref}

\newtheorem{prop}{Proposition}
\newtheorem{thm}{Theorem}
\theoremstyle{remark}
\newtheorem{cor}{Corollary}
\newtheorem{rmk}{Remark}

\makeatletter
\long\def\@makecaption#1#2{%
  \vskip\abovecaptionskip
  \sbox\@tempboxa{\small\textbf{#1:} #2}%
  \ifdim\wd\@tempboxa>\hsize
    {\small\textbf{#1:} #2\par}%
  \else
    \global\@minipagefalse
    \hb@xt@\hsize{\hfil\box\@tempboxa\hfil}%
  \fi
  \vskip\belowcaptionskip}
\makeatother

\title{Testing the Pareto law through a Stein-Laplace characterization}
\author{\parbox{\dimexpr\textwidth-2\tabcolsep\relax}{\centering
Deepesh Bhati$^{1}$, Bruno Ebner$^{2}$, Sakshi Khandelwal$^{1}$\\[0.5ex]
{\small $^{1}$Department of Statistics, Central University of Rajasthan, Ajmer, India\\
$^{2}$Institute of Stochastics, Karlsruhe Institute of Technology (KIT), Karlsruhe, Germany.}
}}
\date{}

\begin{document}


\baselineskip13pt

\maketitle

\begin{abstract}
A Stein--Laplace characterization yields two closed-form goodness-of-fit statistics for the Pareto type-I family with unknown shape and known lower endpoint. The signed integral statistic $T_n$ has a standard normal limit after studentization, but explicit non-Pareto alternatives show that it is not omnibus consistent. The integrated squared discrepancy $S_n$ yields a consistent test against every fixed non-Pareto alternative on the specified support with a finite, positive logarithmic moment. We derive the weighted chi-square null limit of $nS_n$ and give conditions for an asymptotically valid parametric-bootstrap and plug-in covariance-operator method. The latter uses Nystr\"om eigenvalue approximation to avoid repeated sample generation and parameter estimation. A large Monte Carlo comparison study shows competitive performance relative to
existing procedures in the literature. In the paired study, the two calibrations have similar empirical
power, with slightly higher small-sample rejection rates for the bootstrap in some settings. In the timing experiment,
covariance-operator calibration is approximately $28$--$37$ times faster at $n=200$. An exploratory threshold-sensitivity analysis of Danish fire-insurance losses illustrates the dependence of Pareto goodness-of-fit conclusions on the specified tail region.
\end{abstract}


\section{Introduction}

Pareto models describe heavy-tailed data through a simple distributional form above a lower endpoint. Introduced by \citet{pareto1896cours} to describe wealth, they are also used in economics \citep{simon1955class}, finance \citep{rydberg2000realistic}, insurance \citep{embrechts1997modelling}, and environmental science \citep{katz2002statistics}. For a Pareto type-I random variable $X\sim\mathcal P(\alpha)$, where $\alpha>0$ is a shape parameter, the distribution function and probability density function are
\begin{equation}
\label{Pa}
F(x;\alpha) :=
    \begin{cases}
      1-\dfrac{1}{x^\alpha}, & x \geq 1\\
      0, & \text{otherwise}
    \end{cases},
\qquad
f(x;\alpha) :=
    \begin{cases}
     \frac{\alpha}{x^{\alpha+1}}, & x\geq1\\
      0, & \text{otherwise}
    \end{cases},
\quad \alpha>0.
\end{equation}

\noindent Let $X_1,\ldots, X_n$ be independent and identically distributed (iid) random variables with an unknown distribution function $F$ supported on $[1,\infty)$. We test the composite null hypothesis
\[H_0:F\in\{\mathcal P(\alpha):\alpha>0\}\]
for some unknown shape parameter $\alpha>0$, against the alternative
\[H_1:F\notin\{\mathcal P(\alpha):\alpha>0\},\]
where the decision is based on the data $X_1,\ldots, X_n$. For a different positive known lower threshold $x_m>0$, the observations should first be divided by $x_m$.

\noindent Existing Pareto goodness-of-fit procedures use empirical distribution functions \citep{choulakian2001goodness}, probability plots \citep{castillo2012extreme}, or empirical characteristic functions \citep{ndwandwe2023new}. Characterization-based tests begin with an identity that distinguishes the Pareto family from other distributions; see \citet{Ngatchou,bhati2025new}.

\noindent Stein's method, originally developed for normal approximation \citep{stein1972bound}, characterizes probability distributions through operator identities. Using this approach, \citet{betsch2021fixed} derived fixed-point characterizations of continuous univariate distributions. For the Pareto type-I distribution, \citet{avhad2026goodness} established a Stein-type characterization and proposed goodness-of-fit tests based on integral-type and Cram\'er--von Mises-type statistics for both complete and right-censored samples. More recently, \citet{bhati2026laplace} combined the Laplace transform with the Stein-based fixed-point characterization of \citet{betsch2021fixed}. They obtained a characterization of the Pareto type-I distribution and proposed supremum, integral, and $L^2$-type test statistics based on its empirical counterpart.

In the present paper, we use exponential test functions within the Stein density approach to obtain a Laplace-transform like characterization of the Pareto type-I distribution and construct two test statistics. Integrating the signed empirical discrepancy yields $T_n$, a simple moment-type diagnostic, whereas integrating its square yields $S_n$. Squaring before integration prevents cancellation between discrepancies at different transform arguments. This distinction has consequences for consistency of the procedures. The test based on $T_n$ is blind against certain non-Pareto alternatives, whereas the test based on $S_n$ is consistent against every fixed non-Pareto alternative satisfying the logarithmic moment condition stated below.

A studentized version of $T_n$ has a standard normal limiting null distribution, permitting direct asymptotic calibration without resampling. In contrast, the limiting null distribution of $nS_n$ is an infinite weighted sum of independent $\chi^2_1$ random variables. The weights are the eigenvalues of the covariance operator of the limiting Gaussian process and depend on the unknown shape parameter $\alpha_0$. We consider two methods for obtaining critical values. The first uses a parametric bootstrap, with the shape parameter re-estimated in each bootstrap sample. The second replaces $\alpha_0$ by $\widehat{\alpha}_n$ in the covariance operator, numerically approximates its eigenvalues, and uses the resulting weighted chi-square distribution to obtain approximate critical values. We compare the computational cost and finite-sample performance of these approaches in a simulation study. In the settings considered, the plug-in procedure is computationally faster, particularly for larger samples, whereas the bootstrap generally yields slightly higher empirical power for small samples.

\noindent The remainder of the article is organized as follows. Section~\ref{Fusion} presents the characterization and constructs
the test statistics. Section~\ref{Asymptotics} derives their null
limits, identifies the alternatives detected consistently by $T_n$,
and establishes consistency and calibration results for $S_n$.
Section~\ref{MCstudy} compares empirical size and power with those
of existing procedures. Section~\ref{sec:Sn_comparison} compares
the two calibrations in terms of rejection rates and computation
time. Section~\ref{datanalysis} examines threshold sensitivity
using Danish fire-insurance losses. Section~\ref{conclusions}
discusses the implications and directions for further research.

\section{New goodness-of-fit tests for the Pareto law} \label{Fusion}
\subsection{A Stein--Laplace characterization}
We specialize the density approach to Stein's method \citep{ley2013density,ley2017comparison} to the Pareto family \citep[Section~3]{betsch2021fixed}. Exponential test functions yield the following characterization.

\begin{thm}\label{thm:char}
 Let $\alpha>0$ and let $X$ be a positive random variable. Then $X$ has the Pareto type-I distribution $\mathcal P(\alpha)$ if and only if, for every $t>0$,
 \begin{equation}\label{eq:char}
 \alpha\exp(-t)-\mathbb{E}\left[\left(t+\frac{\alpha+1}{X}\right)\exp(-tX)\right]=0.
 \end{equation}
\end{thm}
\begin{proof}
Integration by parts proves the forward implication. For $X\sim\mathcal P(\alpha)$,
\[
\mathbb{E}\!\left[\left(t+\frac{\alpha+1}{X}\right)e^{-tX}\right]
=-\int_1^\infty\frac{\mathrm{d}}{\mathrm{d}x}\{\alpha e^{-tx}x^{-\alpha-1}\}\,\mathrm{d}x
=\alpha e^{-t}.
\]
For the converse, we identify the Laplace transform. Let $\varphi(t):=\mathbb{E}[e^{-tX}]$ and suppose \eqref{eq:char} holds. The nonnegative terms inside the expectation in \eqref{eq:char} are integrable. Tonelli's theorem therefore gives
\[
t\varphi(t)+(\alpha+1)\int_t^\infty\varphi(s)\,\mathrm{d}s=\alpha e^{-t}.
\]
The Laplace transform $\varphi$ is differentiable for $t>0$, since the function $x\mapsto xe^{-tx}$ is bounded for each such $t$. Differentiation gives $t\varphi'(t)-\alpha\varphi(t)=-\alpha e^{-t}$. Multiplying by $t^{-\alpha-1}$ and integrating to infinity, using $0\leq\varphi(t)\leq1$, gives
\[
\varphi(t)=\alpha t^\alpha\int_t^\infty u^{-\alpha-1}e^{-u}\,\mathrm{d}u
=\int_1^\infty e^{-tx}\alpha x^{-\alpha-1}\,\mathrm{d}x.
\]
The uniqueness theorem for Laplace transforms \citep[Chapter~II]{widder1941laplace} identifies the law of $X$.
\end{proof}

The boundary term in Theorem~\ref{thm:char} comes from the lower endpoint of the Pareto support. To see its relation to the density approach, write $p_\alpha(x):=\alpha x^{-\alpha-1}\mathbf{1}_{\{x\geq1\}}$, where $\mathbf{1}_A$ denotes the indicator of a set or event $A$. On the interior of the support, the standard density-Stein operator \citep[Section~3]{betsch2021fixed} acts on a test function $g$ as
\[
(\mathcal A_\alpha g)(x):=g'(x)-\frac{\alpha+1}{x}g(x).
\]
For suitable $g$, integration by parts gives
\[
\int_1^\infty\{g'(x)p_\alpha(x)+g(x)p_\alpha'(x)\}\,\mathrm{d}x=-\alpha g(1).
\]
Choosing the characterizing function class $\{g_t(x)=e^{-tx}:\;t>0\}$ yields \eqref{eq:char}.

\subsection{Two ways to aggregate the discrepancy}

Let $X_1, X_2,\ldots, X_n$ be iid with distribution function $F$ on $[1,\infty)$. We estimate the shape parameter $\alpha$ by the Pareto maximum likelihood estimator (MLE) $\widehat{\alpha}_n:=n/\sum_{j=1}^n\log X_j$, for details see \citet[Eq. (20.42)]{johnson1994continuous}.
\begin{rmk}
     For alternatives that may be identical to $1$ almost surely, we set $\widehat\alpha_n=1$ on the event $\sum_j\log X_j=0$ and define the statistics using this convention. If $\mathbb{E}_F[\log X_1]>0$, then $F(\{1\})<1$ and this event occurs only finitely often almost surely. This does not affect any limit below. Under the null hypothesis, it has probability zero.
\end{rmk}
\noindent We replace the expectation in \eqref{eq:char} by the sample average and substitute $\widehat\alpha_n$ for $\alpha$, which leads to the empirical discrepancy
\begin{equation}\label{eq:Dn}
D_n(t):=\widehat{\alpha}_n e^{-t}
-\frac{1}{n}\sum_{j=1}^{n}
\bigl(t+(\widehat{\alpha}_n+1)/X_j\bigr)e^{-tX_j},\qquad t>0.
\end{equation}
Integrating gives
\begin{align}\label{Tn}
T_n&:=\int_0^\infty D_n(t)\,\mathrm{d}t=\widehat\alpha_n-(\widehat\alpha_n+2)\frac1n\sum_{j=1}^nX_j^{-2}.
\end{align}
This statistic reduces the transform identity to a single moment restriction. To retain discrepancies across all values of $t$, we square before integrating and obtain
\begin{eqnarray}\label{Sn}
S_n:=\int_0^\infty\left[D_n(t)\right]^2\,\mathrm{d}t&=&\frac{\widehat\alpha_n^2}{2}
-\frac{2\widehat\alpha_n}{n}\sum_{j=1}^n\left(\frac1{(1+X_j)^2}
+\frac{\widehat\alpha_n+1}{X_j(1+X_j)}\right)\nonumber\\
&&\quad+\frac1{n^2}\sum_{j,k=1}^n\left(\frac2{(X_j+X_k)^3}
+\frac{(\widehat\alpha_n+1)(\widehat\alpha_n+2)}{X_jX_k(X_j+X_k)}\right).
\end{eqnarray}
Both expressions can be evaluated efficiently without numerical integration. The statistic $T_n$ can depart from zero in either direction, whereas $S_n$ is nonnegative and large values indicate a departure from the characterization.

\begin{rmk}
To compare our construction with that of \citet{bhati2026laplace},
note that their empirical discrepancy can be written as
\[R_n(t):=\frac{1}{n}\sum_{j=1}^{n}
\left[
(\widehat{\alpha}_n+1)
\frac{e^{-t}-e^{-tX_j}}{tX_j}
-e^{-tX_j}
\right],\quad t>0.\]
The two discrepancies satisfy the exact identity
$D_n(t)=tR_n(t)+r_n e^{-t}$, with
$r_n:=\widehat{\alpha}_n
-\frac{\widehat{\alpha}_n+1}{n}
\sum_{j=1}^{n}X_j^{-1}$. The maximum likelihood estimator does not make $r_n$ identically
zero, and the factor $t$ changes the weighting across transform
arguments. In particular, their signed integral is
$\int_0^\infty R_n(t)\,\mathrm{d}t
=\frac{1}{n}\sum_{j=1}^{n}
((\widehat{\alpha}_n+1)\log X_j-1)/X_j$,
whereas ours is given by \eqref{Tn}.
Likewise, their integrated squared discrepancy is
$\int_0^\infty R_n(t)^2\,\mathrm{d}t$,
rather than \eqref{Sn}.
\end{rmk}

\section{Asymptotics} \label{Asymptotics}
For the limit null results, assume that $X_1,\ldots,X_n$ are iid from $\mathcal P(\alpha_0)$, where $\alpha_0>0$. The fixed-alternative results state their assumptions separately. Here, a subscript $F$ on an expectation, variance, covariance, or probability specifies the sampling distribution. Under the null, $\log X_1$ has an exponential distribution with rate $\alpha_0$. The standard delta method for the reciprocal sample mean \citep[Chapter~3]{Vaart_1998} gives the asymptotic linear Bahadur representation
\begin{equation}\label{eq:Bah}
\sqrt{n}\left(\widehat\alpha_n-\alpha_0\right)  = -\frac{\alpha_0^2}{\sqrt{n}}\sum_{i=1}^n\Bigl(\log X_i-\frac{1}{\alpha_0}\Bigr)+o_p(1).
\end{equation}
No finite positive moment of $X_1$ is required for this expansion.
\subsection{Null asymptotics of \texorpdfstring{$T_n$}{Tn}}
The statistic $T_n$ depends on the sample only through the means of $\log X_j$ and $X_j^{-2}$. This reduces its asymptotic analysis to a bivariate central limit theorem and the delta method.

\begin{thm} Under the stated assumptions we have for $n\to\infty$
\begin{equation*}
    \sqrt{n}T_n\stackrel{\mathcal{D}}{\rightarrow} \mathrm{N}(0,\sigma^2),
\end{equation*}
where $\sigma^2:=16\alpha_0/\{(\alpha_0+4)(\alpha_0+2)^2\}$.
\end{thm}
\begin{proof}
Under the null hypothesis,
\[
\mathbb{E}[\log X_1]=\frac1{\alpha_0},\quad
\mathbb{E}[X_1^{-2}]=\frac{\alpha_0}{\alpha_0+2},\]
and
\[
\mathbb{V}[\log X_1]=\frac1{\alpha_0^2},\quad
\mathbb{V}[X_1^{-2}]=\frac{4\alpha_0}{(\alpha_0+4)(\alpha_0+2)^2},
\]
as well as $\mathrm{Cov}[\log X_1,X_1^{-2}]=-2/(\alpha_0+2)^2$.

\noindent The function $h$ defined by $h(x,y):=(1-y)/x-2y$ maps the two sample means to $T_n$. The multivariate delta method \citep[Chapter~3]{Vaart_1998} therefore gives the centered influence variable
\[
\xi:=-\frac{2\alpha_0^2}{\alpha_0+2}\left(\log X_1-\frac1{\alpha_0}\right)
-(\alpha_0+2)\left(X_1^{-2}-\frac{\alpha_0}{\alpha_0+2}\right).
\]
Substitution of these moments gives
$\mathbb{V}[\xi]=16\alpha_0/\{(\alpha_0+4)(\alpha_0+2)^2\}$ and hence the asserted limit.
\end{proof}
\noindent The limiting variance depends on the unknown shape parameter. Substituting its consistent estimate gives a parameter-free normal limit and hence an asymptotic calibration for $T_n$.
\begin{cor}\label{cor:Tn}
    Under the assumptions of the preceding theorem, let $\widehat\alpha_n=n/\sum_{j=1}^n\log X_j$ be the maximum likelihood estimator used in $T_n$. Then
    \begin{equation*}
        \frac{(\widehat{\alpha}_n+2)}4\sqrt{n\frac{\widehat{\alpha}_n+4}{\widehat{\alpha}_n}}T_n\stackrel{\mathcal{D}}{\rightarrow} \mathrm{N}(0,1),
    \end{equation*}
    as $n\rightarrow\infty$.
\end{cor}
\noindent The corollary follows from Slutsky's theorem. A different positive consistent estimator may be used in the studentizing factor, provided that $T_n$ itself retains the MLE. The simulation study also includes a bootstrap procedure for comparability with the other statistics.

\subsection{Fixed alternatives and limits of consistency}
The same two sample means determine the behavior of $T_n$ under a fixed alternative. Their population counterparts yield a limit $\Delta_F$, where a nonzero value leads to consistency, whereas cancellation can make the statistic insensitive to a non-Pareto distribution. The following theorem gives both the almost-sure limit and the fluctuations around it.

\begin{thm}\label{thm:Tn_alt}
Let $X_1,X_2,\ldots$ be iid from a distribution $F$ on $[1,\infty)$ such that $0<\mu_L:=\mathbb{E}_F[\log X_1]<\infty$, $\mathbb{E}_F[(\log X_1)^2]<\infty$. Define $\mu_Y:=\mathbb{E}_F[X_1^{-2}]$, and $\Delta_F:=\mu_L^{-1}-(\mu_L^{-1}+2)\mu_Y$. Then
\[
\widehat\alpha_n \xrightarrow{\mathrm{a.s.}} \mu_L^{-1}
\qquad\text{and}\qquad
T_n \xrightarrow{\mathrm{a.s.}} \Delta_F.
\]
Moreover,
\[
\sqrt n\,(T_n-\Delta_F)\xrightarrow{\mathcal D} \mathrm{N}(0,\tau_F^2),
\]
where
\[
\tau_F^2
=
(1-\mu_Y)^2\mu_L^{-4}\mathbb{V}_F[\log X_1]
+(\mu_L^{-1}+2)^2\mathbb{V}_F[X_1^{-2}]
+2(1-\mu_Y)\mu_L^{-2}(\mu_L^{-1}+2)\mathrm{Cov}_F[\log X_1,X_1^{-2}].
\]
\end{thm}

\begin{proof}
Let $\overline L_n:=n^{-1}\sum_j\log X_j$ and $\overline Y_n:=n^{-1}\sum_jX_j^{-2}$. Except on the event $\sum_j\log X_j=0$, which occurs only finitely often almost surely under the stated assumptions, $\widehat\alpha_n=1/\overline L_n$ and $T_n=h(\overline L_n,\overline Y_n)$. The strong law and the positivity of $\mu_L$ give both almost-sure limits.

\noindent For the normal limit, boundedness of $X_1^{-2}$ and the finite second moment of $\log X_1$ allow us to apply the multivariate central limit theorem to $(\log X_1,X_1^{-2})^\top$. The delta method \citep[Chapter~3]{Vaart_1998} applied at $(\mu_L,\mu_Y)^\top$ provides $\tau_F^2:=\boldsymbol{\nabla}h(\mu_L,\mu_Y)^\top \boldsymbol{\Sigma}_F\,\boldsymbol{\nabla}h(\mu_L,\mu_Y)$ with $h(\ell,y)=(1-y)/\ell-2y$ and covariance matrix
\[
\boldsymbol{\Sigma}_F:=
\begin{pmatrix}
\mathbb{V}_F[\log X_1] & \mathrm{Cov}_F[\log X_1,X_1^{-2}]\\[1mm]
\mathrm{Cov}_F[\log X_1,X_1^{-2}] & \mathbb{V}_F[X_1^{-2}]
\end{pmatrix}.
\]
The asserted limit follows by direct calculation using the gradient $\boldsymbol{\nabla}h(\ell,y)=\left(-(1-y)/\ell^2,-1/\ell-2\right)^\top$.
Substitution yields the displayed variance formula, with a degenerate normal law understood when $\tau_F^2=0$.
\end{proof}

\noindent Let $\Phi$ denote the standard normal distribution function and let $z_p:=\Phi^{-1}(p)$ be its $p$-quantile for $0<p<1$. Fix a significance level $0<\beta<1$. The next corollary identifies the alternatives detected consistently by the two-sided normal test.
\begin{cor}\label{cor:Tn_consistency}
Assume the conditions of Theorem~\ref{thm:Tn_alt}. If $\Delta_F\neq 0$, then $|\,\sqrt n\,T_n\,|\xrightarrow{\mathbb P}\infty$. Define the studentized statistic
\[
W_n:=\frac{\widehat\alpha_n+2}{4}\sqrt{\frac{n(\widehat\alpha_n+4)}{\widehat\alpha_n}}\,T_n.
\]
The asymptotic normal test based on $W_n$ is therefore consistent against every fixed alternative $F$ such that $\Delta_F\neq 0$, that is,
\[
\mathbb{P}_F\left[|W_n|>z_{1-\beta/2}\right]\to 1
\qquad\text{as } n\to\infty.
\]
\end{cor}

\begin{proof}
By Theorem~\ref{thm:Tn_alt}, $T_n\to\Delta_F\ne0$ in probability. The factor multiplying $\sqrt n\,T_n$ in $W_n$ converges to the finite positive number
$(\mu_L^{-1}+2)\sqrt{1+4\mu_L}/4$.
Thus $|\sqrt n\,T_n|\to\infty$ and $|W_n|\to\infty$ in probability, which proves the claim.
\end{proof}
The condition $\Delta_F\neq0$ cannot be omitted. The following proposition constructs continuous non-Pareto alternatives with $\Delta_F=0$ against which the two-sided normal test is not consistent. To establish this failure, we compute the limiting rejection probability, since the almost-sure convergence of $T_n$ to zero alone does not determine the test's power.

\begin{prop}\label{prop:continuous_invisible}
Let $\lambda>0$ and $b>0$, and define
\[
q_{\lambda,b}
:=
\frac{2\lambda b-(\lambda+2)(1-e^{-2b})}{2\lambda b(1-e^{-2b})}.
\]
Assume that $q_{\lambda,b}\in(0,1)$. Let $F_{\lambda,b}$ be the distribution on $[1,\infty)$ with density
\[
f_{\lambda,b}(x)
:=
\lambda x^{-(\lambda+1)}
\Bigl[(1-q_{\lambda,b})+q_{\lambda,b}e^{\lambda b}\mathbf{1}_{\{x\ge e^b\}}\Bigr],
\qquad x\ge 1.
\]
Then $F_{\lambda,b}$ is absolutely continuous, has support $[1,\infty)$, is not a Pareto distribution, and satisfies
\[
\Delta_{F_{\lambda,b}}=0.
\]
Consequently,
\[
T_n \xrightarrow{\mathrm{a.s.}} 0
\qquad\text{under }F_{\lambda,b},
\]
so the test based on $T_n$ is not consistent against $F_{\lambda,b}$.
\end{prop}

\begin{proof}
Let $E\sim\mathrm{Exp}(\lambda)$ and $B\sim\mathrm{Bernoulli}(q_{\lambda,b})$ be independent. The random variable $X:=e^{E+bB}$ has the displayed density. The multiplicative factor in this density changes at $e^b$, so it is not a Pareto density. Writing $q:=q_{\lambda,b}$, the moments satisfy
\[
\mu_L=\frac1\lambda+qb,\qquad
\mu_Y=\frac{\lambda}{\lambda+2}\{1-q(1-e^{-2b})\}
=\frac1{1+2\mu_L},
\]
where the last identity follows by substitution of $q$. Thus $\Delta_F=0$.

\noindent Write $\sigma^2(a):=16a/\{(a+4)(a+2)^2\}$ for the null variance at shape $a$. Since $\log X=E+bB$ has a finite second moment, Theorem~\ref{thm:Tn_alt} gives $T_n\to0$ almost surely and
\[
W_n\xrightarrow{\mathcal D}
\mathrm{N}\!\left(0,\frac{\tau_F^2}{\sigma^2(\mu_L^{-1})}\right).
\]
To see that $0<\tau_F^2<\infty$, consider the influence function $\psi_F$, defined by
\[
\psi_F(y):=-(1-\mu_Y)\mu_L^{-2}(y-\mu_L)
-(\mu_L^{-1}+2)(e^{-2y}-\mu_Y).
\]
This function is nonconstant on $(0,\infty)$, where $E+bB$ has a positive density, and has a finite second moment. Consequently,
\[
\mathbb{P}_F[|W_n|>z_{1-\beta/2}]
\longrightarrow 2\left[1-\Phi\!\left(
\frac{z_{1-\beta/2}\sigma(\mu_L^{-1})}{\tau_F}\right)\right]<1.
\]
The limiting rejection probability is strictly less than one, so the test is not consistent against this alternative.
\end{proof}
\subsection{Null limit of \texorpdfstring{$S_n$}{Sn}}\label{sec:Sn_null}
The null limit of $nS_n$ follows from a central limit theorem for the process inside its squared integral. To express the statistic as a squared norm, we work in the separable Hilbert space $\mathbb{H}:=L^2([0,\infty),\mathcal{B},\mathrm{d}t)$ of equivalence classes of measurable functions $f:[0,\infty) \rightarrow \mathbb{R}$
satisfying $\int_0^\infty |f(t)|^2\,\mathrm{d}t < \infty$, where $\mathcal{B}$ denotes the Borel sigma-field on $[0,\infty)$. Define its scalar product and norm by
\begin{equation*}
\langle f,g \rangle_{\mathbb{H}} := \int_0^\infty f(t)g(t)\,\mathrm{d}t, \quad  \lVert f \rVert_{\mathbb{H}} := \langle f,f \rangle_{\mathbb{H}}^{1/2}, \quad f,g \in \mathbb{H},
\end{equation*}
respectively. The corresponding empirical process is
\begin{equation*}
Z_n(t):=\sqrt{n}\widehat{\alpha}_n\exp(-t)-\frac1{\sqrt{n}}\sum_{j=1}^n\left(t+\frac{\widehat{\alpha}_n+1}{X_j}\right)\exp(-tX_j),\quad t>0,
\end{equation*}
so that $nS_n=\int_0^\infty Z_n^2(t)\,\mathrm{d}t$. Separating the parameter-estimation term gives
\begin{eqnarray*}
Z_n(t)&=&\sqrt{n}\alpha_0\exp(-t)-\frac1{\sqrt{n}}\sum_{j=1}^n\left(t+\frac{\alpha_0+1}{X_j}\right)\exp(-tX_j)\\&&+\sqrt{n}(\widehat{\alpha}_n-\alpha_0)\left(\exp(-t)-\frac1n\sum_{j=1}^n\frac{\exp(-tX_j)}{X_j}\right),\quad t>0.
\end{eqnarray*}
The Bahadur representation \eqref{eq:Bah} leads to the approximation process
\begin{eqnarray*}
\widehat{Z}_n(t)&:=&\sqrt{n}\alpha_0\exp(-t)-\frac1{\sqrt{n}}\sum_{j=1}^n\left(t+\frac{\alpha_0+1}{X_j}\right)\exp(-tX_j)\\&&- \frac{\alpha_0^2}{\sqrt{n}}\sum_{i=1}^n\Bigl(\log X_i-\frac{1}{\alpha_0}\Bigr)\left(\exp(-t)-\Psi_{\alpha_0}(t)\right),\quad t>0,
\end{eqnarray*}
where $\Psi_{\alpha_0}(t):=\alpha_0t^{\alpha_0+1}\Gamma(-\alpha_0-1,t)$ and the upper incomplete gamma function is defined by $\Gamma(s,t):=\int_t^\infty u^{s-1}\exp(-u)\,\mathrm{d}u$.
The next theorem identifies the Gaussian limit.
\begin{thm}\label{thm:AsyDist}
Under the standing assumptions, there exists a centered Gaussian random element $Z$ in $\mathbb{H}$ with covariance kernel
\begin{eqnarray*}
K_{\alpha_0}(s,t)&:=&\alpha_0 st\,(s+t)^{\alpha_0}\Gamma(-\alpha_0,s+t)
 +\alpha_0(\alpha_0+1)(s+t)^{\alpha_0+2}\Gamma(-\alpha_0-1,s+t)\\
&&+\alpha_0(\alpha_0+1)^2(s+t)^{\alpha_0+2}\Gamma(-\alpha_0-2,s+t)
 -\alpha_0^2e^{-s-t}-\alpha_0^2C_{\alpha_0}(s)C_{\alpha_0}(t),
\end{eqnarray*}
where $C_{\alpha_0}(t):=e^{-t}-\Psi_{\alpha_0}(t)$, $s,t>0$, such that
\[
nS_n\xrightarrow{\mathcal{D}}\lVert Z\rVert_{\mathbb{H}}^2,\qquad \mbox{as }n\to\infty.
\]
\end{thm}

\begin{proof}
We first show that the approximating process is a sum of centered, square-integrable random elements and that its difference from $Z_n$ is negligible. We then calculate the covariance kernel.
Throughout the proof let $\mathbb{E}_{\alpha_0}$ denote expectation under $\mathcal P(\alpha_0)$, and define $Y_x(t):=e^{-tx}/x$,
and
\[
h_{\alpha_0}(x,t):=\alpha_0 e^{-t}-\left(t+\frac{\alpha_0+1}{x}\right)e^{-tx}
-\alpha_0^2\left(\log x-\frac{1}{\alpha_0}\right)C_{\alpha_0}(t).
\]
The characterization gives $\mathbb{E}_{\alpha_0}[h_{\alpha_0}(X,t)]=0$. Since $x\ge1$,
$\|Y_x\|_{\mathbb H}^2=1/(2x^3)\le1/2$,
$0\le C_{\alpha_0}(t)\le e^{-t}$, and
\[
|h_{\alpha_0}(x,t)|\le c_{\alpha_0}(1+\log x)(1+t)e^{-t},
\]
where $c_{\alpha_0}$ is a finite constant depending only on $\alpha_0$. Thus $h_{\alpha_0}(X,\cdot)$ is square-integrable in $\mathbb H$, because
$\log X\sim\mathrm{Exp}(\alpha_0)$. The Hilbert-space law of large numbers
and central limit theorem \citep[Chapter~17]{H:2024} therefore apply.
In particular, with $H_j:=h_{\alpha_0}(X_j,\cdot)$ and
$\overline Y_n:=n^{-1}\sum_jY_{X_j}$, we have
$\|\overline Y_n-\mathbb{E}_{\alpha_0}[Y_X]\|_{\mathbb H}\to0$ almost surely and
$n^{-1/2}\sum_jH_j\xrightarrow{\mathcal{D}} Z$ in $\mathbb H$.
If $R_n=o_{\mathbb P}(1)$ is the remainder in \eqref{eq:Bah}, then
\[
Z_n-\frac1{\sqrt n}\sum_{j=1}^nH_j
=R_nC_{\alpha_0}-\sqrt n(\widehat\alpha_n-\alpha_0)
 (\overline Y_n-\mathbb{E}_{\alpha_0}[Y_X])=o_{\mathbb P}(1)
\quad\text{in }\mathbb H.
\]
Slutsky's theorem and the continuous mapping theorem give the asserted
limit of $nS_n=\|Z_n\|_{\mathbb H}^2$.  It remains to calculate the kernel. Set $B_t(x):=\{t+(\alpha_0+1)/x\}e^{-tx}$, suppressing its dependence on $\alpha_0$ in the notation.
Differentiating the identity $\mathbb{E}_{\alpha_0}[B_t(X)]=\alpha_0 e^{-t}$ with respect
to $\alpha_0$ is justified locally by the exponential law of $\log X$ and gives
\[
\mathbb{E}_{\alpha_0}\left[\left(\log X-\frac1{\alpha_0}\right)B_t(X)\right]=-C_{\alpha_0}(t).
\]
This differentiation includes the dependence on $\alpha_0$ of both $B_t$ and the Pareto law.
Consequently, after some more calculations, the kernel is as stated in the theorem.
\end{proof}

\subsection{Consistency of \texorpdfstring{$S_n$}{Sn}}\label{sec:Sn_consistency}
Under the logarithmic moment condition below, $S_n$ converges to the squared norm of the population discrepancy. The characterization in Theorem~\ref{thm:char} makes this limit positive whenever the distribution is not Pareto.

\begin{thm}\label{thm:Sn_consistency}
Let $X_1,X_2,\ldots$ be iid from a distribution $F$ on $[1,\infty)$ such that $0<\mu_L=\mathbb{E}_F[\log X_1]<\infty$, and set $\alpha_F=1/\mu_L$. Define
\begin{equation}\label{eq:population_discrepancy}
D_F(t):=\alpha_F e^{-t}-\mathbb{E}_F\left[\left(t+\frac{\alpha_F+1}{X_1}\right)e^{-tX_1}\right],
\qquad t>0,
\end{equation}
and $S(F):=\int_0^\infty D_F^2(t)\,\mathrm{d}t$. Then
\begin{equation}\label{eq:Sn_population_limit}
S_n\xrightarrow{\mathrm{a.s.}}S(F).
\end{equation}
Moreover, $S(F)=0$ if and only if $F=\mathcal P(\alpha_F)$. Hence, for every fixed non-Pareto alternative satisfying the stated moment condition, $S(F)>0$ and $nS_n\to\infty$ almost surely. Consequently, every right-tailed test based on $nS_n$ with critical values bounded in probability is consistent.
\end{thm}

\begin{proof}
We prove convergence of the empirical discrepancy in $\mathbb H$ and then identify the distributions for which its limit vanishes. The reciprocal-mean formula for $\widehat\alpha_n$ holds eventually almost surely under the convention above, so the strong law of large numbers gives $\widehat\alpha_n\to\alpha_F$ almost surely.

\noindent For fixed $a>0$, define $g_a(x,t):=\{t+(a+1)/x\}e^{-tx}$. Since $x\geq 1$, we have $\int_0^\infty g_a^2(x,t)\,\mathrm{d}t\leq c_a<\infty$ for a constant $c_a$ depending only on $a$. Thus $g_{\alpha_F}(X_1,\cdot)$ is square-integrable in $\mathbb H$, and the Hilbert-space strong law of large numbers \citep[Theorem 17.15]{H:2024} yields
\[
\Bigl\lVert n^{-1}\sum_{j=1}^n g_{\alpha_F}(X_j,\cdot)-\mathbb{E}_F[g_{\alpha_F}(X_1,\cdot)]\Bigr\rVert_{\mathbb H}\to 0
\quad\text{almost surely}.
\]
To account for parameter estimation, use the bound $\lVert g_a(x,\cdot)-g_b(x,\cdot)\rVert_{\mathbb H}\leq |a-b|/\sqrt{2}$, which holds uniformly in $x\geq 1$. Replacing $\alpha_F$ by $\widehat\alpha_n$ is therefore asymptotically negligible in $\mathbb H$. The empirical discrepancy defining $S_n$ consequently converges almost surely in $\mathbb H$ to $D_F$, and continuity of the squared norm gives \eqref{eq:Sn_population_limit}.

\noindent If $S(F)=0$, then $D_F(t)=0$ for Lebesgue-a.e. $t>0$. The function $D_F$ is continuous in $t$, so $D_F(t)=0$ for every $t>0$. Theorem~\ref{thm:char}, applied with $\alpha=\alpha_F$, then implies $F=\mathcal P(\alpha_F)$. The converse is immediate from Theorem~\ref{thm:char}. If $F$ is not Pareto, $S(F)>0$ and $nS_n\to\infty$ almost surely.
\end{proof}

\subsection{Covariance-operator calibration}\label{sec:Sn_covariance}
The null limit in Theorem~\ref{thm:AsyDist} depends on the unknown parameter $\alpha_0>0$. We therefore use two different calibrations for $S_n$. Here, we approximate the covariance eigenvalues and substitute $\widehat\alpha_n$ for $\alpha_0$, and Section~\ref{sec:bootstrap_calibration} treats the parametric bootstrap. Let $\mathcal K_{\alpha_0}:\mathbb H\to\mathbb H$ denote the covariance operator defined by
\[
(\mathcal K_{\alpha_0}f)(s)
:=
\int_0^\infty K_{\alpha_0}(s,t)f(t)\,\mathrm{d}t.
\]
The moment bound in the proof of Theorem~\ref{thm:AsyDist} shows that this positive covariance operator is trace class, since
\[
\operatorname{tr}(\mathcal K_{\alpha_0})
=\mathbb{E}\lVert Z\rVert_{\mathbb H}^2
=\mathbb{E}\lVert h_{\alpha_0}(X_1,\cdot)\rVert_{\mathbb H}^2<\infty.
\]
Let $\lambda_1,\lambda_2,\ldots$ denote its nonnegative eigenvalues, ordered nonincreasingly with multiplicity and extended by zeros if necessary. The trace-class property gives $\sum_j\lambda_j<\infty$. For iid standard normal random variables $N_1,N_2,\ldots$, the standard Karhunen--Lo\'eve expansion for Gaussian random elements \citep[Chapter~3]{DW:2014} gives
\[
\lVert Z\rVert_{\mathbb H}^2\stackrel{\mathcal{D}}{=}\sum_{j=1}^\infty \lambda_j N_j^2.
\]
The nonzero eigenvalues correspond to eigenfunctions $f_j$ satisfying the integral equation
\begin{equation}\label{int:eq}
\lambda f(s)=(\mathcal K_{\alpha_0}f)(s),\qquad s > 0.
\end{equation}
We approximate these eigenvalues using the Nystr\"om method (for details see \citet{nystrom1930praktische} and \citet[Chapter~4]{A:1997}), so truncation and quadrature reduce the integral equation to a finite matrix eigenvalue problem. For a different approach to covariance-operator eigenvalues in goodness-of-fit testing, see the Rayleigh--Ritz method of \citet{EJM:2025}.

\noindent Since the integration domain is unbounded, we first truncate the operator to a sufficiently large interval $[0,T]$ and consider
\[
(\mathcal K_{\alpha_0,T}f)(s)
:=
\int_0^T K_{\alpha_0}(s,t)f(t)\,\mathrm{d}t,
\qquad 0\le s\le T.
\]
Next, let $t_1,\dots,t_m\in[0,T]$ be quadrature nodes with positive weights $w_1,\dots,w_m$, so that
\[
\int_0^T g(t)\,\mathrm{d}t
\approx
\sum_{j=1}^m w_j g(t_j).
\]
Then the Nystr\"om approximation of $\mathcal K_{\alpha_0,T}$ is
\[
(\mathcal K_{\alpha_0,T}f)(t_i)
\approx
\sum_{j=1}^m w_j\,K_{\alpha_0}(t_i,t_j)\,f(t_j),
\qquad i=1,\dots,m.
\]
Writing $u_j$ for the discrete counterpart of $f(t_j)$, we replace the continuous eigenvalue problem in \eqref{int:eq} by
\[
\sum_{j=1}^m w_j\,K_{\alpha_0}(t_i,t_j)\,u_j
=
\lambda u_i,
\qquad i=1,\dots,m.
\]
Since the kernel is symmetric, we use the symmetrized system with the column vector $\mathbf{v}:=(\sqrt{w_1}u_1,\ldots,\sqrt{w_m}u_m)^\top$:
\[
\mathbf{B}\mathbf{v}=\lambda\mathbf{v},
\qquad
\mathbf{B}:=(B_{ij})_{i,j=1}^m,
\]
where the entries are defined by $B_{ij}:=\sqrt{w_i}\,K_{\alpha_0}(t_i,t_j)\,\sqrt{w_j}$, $i,j=1,\dots,m$. The eigenvalues of $\mathbf{B}$ approximate those of $\mathcal K_{\alpha_0}$. For the fitted model, we substitute $\widehat\alpha_n$ for $\alpha_0$ and denote the resulting eigenvalues by $\widehat\lambda_{j,n}$ after replacing negative numerical values by zero. Given independent chi-square variables $\chi^2_{1,j}$ with one degree of freedom, we approximate the quantile of $\sum_j\widehat\lambda_{j,n}\chi^2_{1,j}$ by Monte Carlo simulation.

\noindent The next theorem states sufficient conditions for this plug-in method to have the intended asymptotic level. It applies standard weak convergence and quantile convergence \citep[Section~21.1]{Vaart_1998}. The key assumption is $\ell^1$-consistency of the entire ordered eigenvalue sequence, including control of its tail; convergence of finitely many eigenvalues alone is insufficient.

\begin{thm}\label{thm:Nystrom_level}
Assume the null hypothesis $H_0:X_1,\dots,X_n\stackrel{\text{iid}}{\sim}\mathcal P(\alpha_0)$ with $\alpha_0>0$, and suppose
\[
nS_n \xrightarrow{\mathcal D} Q_{\alpha_0}:=\sum_{j=1}^\infty\lambda_j(\alpha_0)\chi_{1,j}^2,
\]
where $(\chi_{1,j}^2)_{j\ge1}$ are iid chi-square random variables with one degree of freedom and $(\lambda_j(\alpha_0))_{j\ge1}$ are the eigenvalues of the covariance operator $\mathcal K_{\alpha_0}$.

\noindent Let $\widehat\alpha_n$ be a consistent estimator of $\alpha_0$, and let $\widehat\lambda_{1,n}\ge \widehat\lambda_{2,n}\ge \cdots \ge 0$
be the eigenvalues obtained from a plug-in Nystr\"om approximation of $\mathcal K_{\widehat\alpha_n}$, extended by zeros. Assume that
\[
\sum_{j=1}^\infty\bigl|\widehat\lambda_{j,n}-\lambda_j(\alpha_0)\bigr|\xrightarrow{\mathbb P}0.
\]
Define
\[
\widehat Q_n:=\sum_{j=1}^\infty\widehat\lambda_{j,n}\chi_{1,j}^2,
\]
conditionally on the data, using auxiliary chi-square variables independent of the observations. Let $\widehat q_{n,1-\beta}$ be its conditional lower $(1-\beta)$-quantile, where $0<\beta<1$.
Write $F_{Q,\alpha_0}$ for the distribution function of $Q_{\alpha_0}$ and set $q:=q_{1-\beta}(\alpha_0):=\inf\{x:F_{Q,\alpha_0}(x)\ge1-\beta\}$. Assume that $F_{Q,\alpha_0}$ is continuous at $q$ and, for every $\varepsilon>0$,
\[
F_{Q,\alpha_0}(q-\varepsilon)<1-\beta<F_{Q,\alpha_0}(q+\varepsilon).
\]
Then
\[
\widehat q_{n,1-\beta}\xrightarrow{\mathbb P} q_{1-\beta}(\alpha_0),
\]
and
\[
\lim_{n\rightarrow\infty}\mathbb{P}_{H_0}\left[nS_n>\widehat q_{n,1-\beta}\right]=
\beta.
\]
\end{thm}

\begin{proof}
The eigenvalue assumption controls the difference between the two weighted chi-square variables under a common coupling. Couple $\widehat Q_n$ and $Q_{\alpha_0}$ with the same independent chi-square variables. Conditionally on the data,
\[
\mathbb{E}\bigl[|\widehat Q_n-Q_{\alpha_0}|\mid X_1,\ldots,X_n\bigr]
\le\sum_{j=1}^\infty|\widehat\lambda_{j,n}-\lambda_j(\alpha_0)|
\xrightarrow{\mathbb P}0.
\]
The bound proves conditional weak convergence in probability. The standard criterion for quantile convergence \citep[Section~21.1]{Vaart_1998}, applied along almost-surely convergent subsequences, then gives the stated convergence of quantiles. Slutsky's theorem and continuity of $F_{Q,\alpha_0}$ at $q$ give the asserted level.
\end{proof}

The quantile condition in Theorem~\ref{thm:Nystrom_level} holds for the present covariance operator. Indeed, $C_{\alpha_0}(t)>0$ for $t>0$, and $h_{\alpha_0}(x,t)$ has a nonzero logarithmic term, so the operator is nonzero. Its summable weighted chi-square law has a continuous distribution function that is strictly increasing on $(0,\infty)$.

The following proposition gives sufficient conditions for
the numerical approximation to satisfy the spectral
convergence assumption of Theorem \ref{thm:Nystrom_level}.

\begin{prop}
\label{prop:numerical-calibration}
Assume $H_0$ and $\widehat{\alpha}_n\to\alpha_0>0$ in probability. Let $(t_{i,n},w_{i,n})_{i=1}^{m_n}$ be the Gauss--Legendre rule on $[0,\tau_n]$, and let
$\widetilde K_n$ be a symmetric numerical approximation to $K_{\widehat{\alpha}_n}$. Define $\epsilon_n:=\max_{i,j\le m_n} 
|\widetilde K_n(t_{i,n},t_{j,n})
-K_{\widehat{\alpha}_n}(t_{i,n},t_{j,n})|$.
Suppose
\begin{equation}
\label{eq:numerical-resolution}
\tau_n\longrightarrow\infty,\qquad
\frac{\tau_n^2}{\sqrt{m_n}}\longrightarrow0,\qquad
m_n\tau_n\epsilon_n\xrightarrow{\mathbb P}0.
\end{equation}
Let $\widehat\lambda_{j,n}$ be the ordered eigenvalues of $\widetilde B_n=
(\sqrt{w_{i,n}w_{j,n}}\,
\widetilde K_n(t_{i,n},t_{j,n}))_{i,j=1}^{m_n}$,
with negative eigenvalues replaced by zero and the sequence extended by zeros. Then $\sum_{j\ge1}|\widehat\lambda_{j,n}
-\lambda_j(\alpha_0)|\to0$ in probability. Consequently, calibration using empirical quantiles from $M_n\to\infty$ conditionally independent weighted chi-square draws has asymptotic level $\beta$.
\end{prop}

\begin{proof}
Throughout the proof, $\|\cdot\|_1$ denotes the trace norm and $\|\cdot\|_{\mathrm{op}}$ the operator norm. Let $G=L^2([0,\infty),e^{-z}\,\mathrm dz)$ and set $g_a(t)(z):=h_a(e^{z/a},t)$, where $h_a$ is the
influence function in the proof of Theorem~4. Then $K_a(s,t)=\langle g_a(s),g_a(t)\rangle_G$. Writing $(v\otimes v)f:=\langle f,v\rangle_Gv$, define
$U_a(t):=g_a(t)\otimes g_a(t)$ and $\mathcal C_a:=\int_0^\infty U_a(t)\,\mathrm dt$.
The operators $\mathcal C_a$ and $\mathcal K_a$ have the same nonzero eigenvalues, as do $\sum_iw_{i,n}U_a(t_{i,n})$ and the exact Nystr\"om matrix $B_n(a)$.

On every compact $A\subset(0,\infty)$, the expression for $h_a$ and its derivative in $t$ give $\|g_a(t)\|_G+\|\partial_tg_a(t)\|_G
\le c_A(1+t)e^{-t}$. Thus $U_a$ is uniformly Lipschitz in trace norm. Moreover,  $\|v\otimes v-w\otimes w\|_1
\le(\|v\|+\|w\|)\|v-w\|$ and dominated convergence show that $a\mapsto\mathcal C_a$ is trace-norm continuous.

Bernstein polynomial approximation of degree $m_n$ approximates $U_a$ uniformly on $[0,\tau_n]$ with trace-norm error at most $c_A\tau_n/\sqrt{m_n}$.
Positivity and polynomial exactness of the quadrature, together with $\sum_iw_{i,n}=\tau_n$, therefore yield
\begin{equation}
\label{eq:numerical-trace-bound}
\sup_{a\in A}
\left\|
\sum_{i=1}^{m_n}w_{i,n}U_a(t_{i,n})-\mathcal C_a
\right\|_1
\le c_A\left\{
\frac{\tau_n^2}{\sqrt{m_n}}
+(1+\tau_n)^2e^{-2\tau_n}
\right\}.
\end{equation}
Combining \eqref{eq:numerical-trace-bound},
\eqref{eq:numerical-resolution}, and plug-in
consistency gives trace-norm convergence to
$\mathcal C_{\alpha_0}$ in probability. Hence the
exact Nystr\"om eigenvalues converge coordinatewise
and their sums converge. Nonnegativity then gives
convergence in $\ell^1$.

For the numerically evaluated matrix,
$\|\widetilde B_n-B_n(\widehat\alpha_n)\|_{\mathrm{op}}
\le\tau_n\epsilon_n$.
The min--max principle bounds the sum of the absolute
eigenvalue differences by $m_n\tau_n\epsilon_n$;
replacing negative eigenvalues by zero cannot
increase their distance from the nonnegative exact
eigenvalues. This proves the asserted spectral
convergence.

Finally, conditionally on the data, the variance
of the empirical distribution function of the
$M_n$ reference draws at any fixed argument is
at most $1/(4M_n)$. Applying Chebyshev's inequality
at points on either side of the limiting quantile,
and using the quantile separation established in
Theorem~6, proves empirical quantile consistency.
The asserted asymptotic level follows.
\end{proof}

\begin{rmk}
With exact kernel evaluation, $\epsilon_n=0$.
For Gauss--Laguerre evaluation, the substitution
$X=e^{z/a}$ expresses the kernel expectations as
integrals against $e^{-z}\,\mathrm dz$ of continuous
functions bounded by $c_A(1+z)^2$ on compact
parameter sets. The positive Gauss--Laguerre
measures converge weakly to this exponential
measure by moment exactness and moment determinacy;
their exact fourth moments ensure uniform
integrability. Consequently, the kernel
approximations converge uniformly on compact
$(a,s,t)$-sets. Writing $K_a^{(L)}$ for the approximation with $L$
Gauss--Laguerre nodes, choose $L_n\to\infty$ so that
$\sup_{a\in[1/n,n],\,s,t\in[0,\tau_n]}
|K_a^{(L_n)}(s,t)-K_a(s,t)|
\le (n m_n\tau_n)^{-1}$.
Such a choice is possible by the compact-uniform convergence.
Since $\mathbb P(\widehat\alpha_n\in[1/n,n])\to1$,
the last condition in \eqref{eq:numerical-resolution} follows.
\end{rmk}

Proposition~\ref{prop:numerical-calibration} concerns an increasing-resolution approximation, while the numerical studies below use fixed-resolution implementations with the settings specified in the respective sections.

\subsection{Parametric-bootstrap calibration}\label{sec:bootstrap_calibration}
The bootstrap generates samples from the fitted Pareto law and repeats parameter estimation for each sample. Its validity follows from stability of the null limits along convergent sequences of shape parameters. The argument applies to both $T_n$ and $S_n$.

\begin{thm}\label{thm:bootstrap_validity}
Let $X_1,\ldots,X_n$ be iid from $\mathcal P(\alpha_0)$, $\alpha_0>0$, and let $\widehat\alpha_n$ be the MLE. For test calibration use exact conditional quantiles, or empirical conditional quantiles from $B_n\to\infty$ independent bootstrap samples. Conditionally on the observations, let $X_1^*,\ldots,X_n^*$ be iid from $\mathcal P(\widehat\alpha_n)$ and compute $T_n^*$ and $S_n^*$ from the bootstrap sample, including re-estimation of $\alpha$. Then the conditional distribution of $nS_n^*$ converges in probability to $Q_{\alpha_0}=\lVert Z\rVert_{\mathbb H}^2$. In addition, the conditional distribution of the studentized version of $T_n^*$ converges in probability to the standard normal law. Thus the parametric bootstrap tests based on $T_n$ and $S_n$ have asymptotic level $\beta$ for every $0<\beta<1$ at continuity points of the corresponding limiting distributions.
\end{thm}

\begin{proof}
We verify the local-parameter convergence condition in the standard
parametric-bootstrap consistency theory \citep[Section~2]{Beran1997};
see also \citet[Chapter~23]{Vaart_1998}. It suffices to establish the null limits along every deterministic parameter sequence converging to $\alpha_0$. At parameter $a>0$, the Pareto family with its fixed
support is a regular one-parameter exponential family, with score
$a^{-1}-\log x$ and Fisher information $a^{-2}$, and its MLE is
asymptotically normal by \eqref{eq:Bah}.

\noindent For any deterministic sequence $a_n$ converging to $a>0$, couple a row of Pareto
observations as $X_{n,j}:=\exp(E_j/a_n)$, where the $E_j$ are iid
$\mathrm{Exp}(1)$. This representation places all rows on the same probability space. The corresponding MLE is $\widetilde a_n:=a_n/\overline E_n$, where $\overline E_n:=n^{-1}\sum_{j=1}^n E_j$, and
\[
\sqrt n(\widetilde a_n-a_n)
=-\frac{a_n}{\sqrt n}\sum_{j=1}^n(E_j-1)
 +\frac{a_n\sqrt n(\overline E_n-1)^2}{\overline E_n}.
\]
The last term is $o_{\mathbb P}(1)$. Let $E$ denote a generic $\mathrm{Exp}(1)$ distributed variable. For $a_n$ in a compact subset of
$(0,\infty)$, the influence functions $h_a$ from the proof of
Theorem~\ref{thm:AsyDist} satisfy, for a finite constant $M$ depending only on that compact set,
\[
|h_{a_n}(e^{E/a_n},t)|\le M(1+E)(1+t)e^{-t},\qquad
\mathbb{E}\bigl[\|h_{a_n}(e^{E/a_n},\cdot)-h_a(e^{E/a},\cdot)\|_{\mathbb H}^2\bigr]\to0.
\]
The second statement follows by dominated convergence. Together with the bound
$\|Y_x\|_{\mathbb H}^2\le1/2$ from the same proof, these bounds verify the negligible
remainder, convergence of covariances and their traces, and Lindeberg conditions for the
Hilbert-space triangular-array central limit theorem
\citep[Chapter~17]{H:2024}. Let $Q_a$ denote the limit in Theorem~\ref{thm:AsyDist} with parameter $a$. It follows that $nS_n\xrightarrow{\mathcal{D}} Q_a$ along every such
sequence. The scalar influence function for $T_n$ is likewise continuous
in $a$ in $L^2$ under this coupling and bounded by $M(1+E)$; therefore
$\sqrt nT_n\xrightarrow{\mathcal{D}}\mathrm{N}(0,\sigma^2(a))$, where $\sigma^2(a)=16a/\{(a+4)(a+2)^2\}$, and studentization converges
to the standard normal law along the same sequences.

\noindent Applying parametric-bootstrap consistency with $a_n$ replaced by
$\widehat\alpha_n\to\alpha_0$ gives both asserted conditional limits.
Conditional quantile consistency and the stipulated quantile conditions
give the level statements. The scalar unstudentized limit also justifies
the equal-tailed bootstrap calibration of $T_n$ used in the simulations.
\end{proof}

\section{Simulation study} \label{MCstudy}
We use sample sizes $n=20,50,100$, $N=10{,}000$ replications, and nominal significance level $\eta=0.05$. For the integrated statistic $S_n$, we report two calibrations: the procedure $S_{n,K}$ uses the plug-in covariance operator described after Theorem~\ref{thm:Nystrom_level}, and the procedure $S_{n,B}$ uses the parametric bootstrap of Theorem~\ref{thm:bootstrap_validity}. For $T_n$, we use a two-sided bootstrap test to facilitate comparison with the competing procedures.

\noindent  The proposed procedures are implemented as follows.
\begin{enumerate}
    \item  Draw a random sample $X_1,\ldots,X_n$ from the distribution under study.
    \item  Estimate the shape parameter by the MLE $\widehat{\alpha}_n=n/\sum_{i=1}^n\log X_i$.
    \item  Evaluate $T_n$ and $S_n$ from \eqref{Tn} and \eqref{Sn}. The statistic $T_n$ is calibrated as a two-sided statistic and $S_n$ as a right-tailed statistic.
    \item  For the bootstrap calibration, generate $B=500$ samples of size $n$ from $\mathcal P(\widehat\alpha_n)$. Re-estimate $\alpha$ in each bootstrap sample and compute $T_n^{*(b)}$ and $S_n^{*(b)}$, $b=1,\ldots,B$. Denote the empirical quantiles of the bootstrap $T_n$ values at levels $\eta/2$ and $1-\eta/2$ by $q^{*T}_{n,\eta/2}$ and $q^{*T}_{n,1-\eta/2}$, respectively. Denote the empirical quantile of the bootstrap $S_n$ values at level $1-\eta$ by $q^{*S}_{n,1-\eta}$.
    \item  For the covariance-operator calibration of $S_n$, compute the eigenvalues of the plug-in Nystr\"om matrix with $\widehat\alpha_n$, $T=12$ and $m=50$, and approximate the $1-\eta$ quantile $q^K_{n,1-\eta}$ of $\sum_j\widehat\lambda_{j,n}\chi^2_{1,j}$ by $20{,}000$ Monte Carlo draws.    \item  Reject with the bootstrap version of $T_n$ if $T_n<q^{*T}_{n,\eta/2}$ or $T_n>q^{*T}_{n,1-\eta/2}$; reject with $S_{n,B}$ if $S_n>q^{*S}_{n,1-\eta}$; reject with $S_{n,K}$ if $nS_n>q^K_{n,1-\eta}$.
\end{enumerate}
The competing procedures defined below use the same simulated samples and the same bootstrap sample size $B=500$. After $N$ replications, we report the percentage of samples in which each test rejects. If $\widehat p$ denotes the rejection proportion, the Monte Carlo standard error of the reported percentage $100\widehat p$ is $100\{\widehat p(1-\widehat p)/N\}^{1/2}$. With $N=10{,}000$, it is about $0.22$ percentage points at $5\%$ and at most $0.50$ percentage points.

The comparison includes tests based on empirical distribution functions, transforms, and powers and minima of the observations. To fix the common notation, let $X_{(1)}\le\cdots\le X_{(n)}$ denote the ordered sample. The empirical distribution function is
\begin{equation*}
	\widehat{F}_n(t) := \frac{1}{n} \sum_{i=1}^{n}\mathbf{1}_{\{X_i\le t\}}.
\end{equation*}

\paragraph{Distribution-function tests.}
The test of \citet{obradovic2015goodness} compares the empirical distribution of the observations with that of the maximum ratio within each pair. The latter is
\begin{equation*}
	M_n(t):=\binom{n}{2}^{-1} \sum_{i=1}^{n-1}\sum_{j=i+1}^{n}\mathbf{1}_{\left\{\max\left(X_i/X_j,X_j/X_i\right)\le t\right\}},\quad t\ge 1.
\end{equation*}
The resulting integrated difference is
\begin{equation*}
	O_n:=\int_{1}^{\infty} \left(M_n(t)-\widehat{F}_n(t)\right)\,\mathrm{d}\widehat{F}_n(t).
\end{equation*}

\noindent The three classical empirical-distribution tests compare the sample with the fitted distribution $F_{\widehat\alpha_n}(\cdot):=F(\cdot,\widehat\alpha_n)$ from \eqref{Pa}. The Kolmogorov--Smirnov statistic is defined by
\[KS_n := \max \bigg(\max_{1 \leq j \leq n} \bigg( \frac{j}{n} - F_{\widehat{\alpha}_n} (X_{(j)}) \bigg) ,  \max_{1 \leq j \leq n} \bigg(F_{\widehat{\alpha}_n} (X_{(j)})-\frac{j-1}{n}\bigg)     \bigg).\]
The Cram{\'e}r--von Mises statistic is
\[ CM_n := \frac{1}{12n} + \sum_{j=1}^n \bigg( F_{\widehat{\alpha}_n} (X_{(j)}) - \frac{2j-1}{2n}\bigg)^2,  \]
whereas the Anderson--Darling statistic is given by
\[ AD_n := -n - \frac{1}{n} \sum_{j=1}^n (2j-1) \bigg( \log(F_{\widehat{\alpha}_n} (X_{(j)})) + \log (1-F_{\widehat{\alpha}_n} (X_{(n+1-j)})) \bigg). \]

\paragraph{Transform tests.}
For the Mellin-transform test of \citet{meintanis2009unified}, take a tuning parameter $a>0$ and define, for $t\ge1$, the auxiliary functions
\[
I_a^{(0)}(t) := (a+\log t)^{-1}, \qquad  I_a^{(1)}(t):= \frac{1-a-\log t}{(a+\log t)^{2}},
\]
and
\[
I_a^{(2)}(t):=	\frac{2-2a+a^2+2(a-1) \log t + \log^2 t}{(a+\log t)^{3}}.
\]
The statistic combines these functions as
\begin{align*}
G_{n,a}&:=\frac1n\sum_{j,k=1}^n\bigg[(\widehat\alpha_n+1)^2I_a^{(0)}(X_jX_k)
+I_a^{(2)}(X_jX_k)\\
&\qquad\qquad\quad+2(\widehat\alpha_n+1)I_a^{(1)}(X_jX_k)\bigg]
+n\widehat\alpha_n^2I_a^{(0)}(1)\\
&\quad-2\widehat\alpha_n\sum_{j=1}^n\big[(\widehat\alpha_n+1)I_a^{(0)}(X_j)+I_a^{(1)}(X_j)\big].
\end{align*}
We use the tuning parameter $a=2$.

\noindent The test of \citet{meintanis2014class} applies the empirical characteristic function to the fitted probability transforms
\[ \widehat U_j := F_{\widehat{\alpha}_n}(X_j)=1-X_{j}^{-\widehat{\alpha}_n} ,\quad j=1,2,\ldots,n.\]
In terms of $\widehat U_1, \ldots, \widehat U_n$, the statistic is
\begin{multline*}
	 M_{n,a} := \frac{1}{n} \sum_{j,k=1}^{n} \frac{2a}{(\widehat U_j-\widehat U_k)^2+a^2} + 2n \left( 2 \tan^{-1} \left(\frac{1}{a}\right) - a \log \left( 1+ \frac{1}{a^2}\right)\right) \\ - 4\sum_{j=1}^{n} \left( \tan^{-1} \left(\frac{\widehat U_j}{a}\right) + \tan^{-1} \left(\frac{1-\widehat U_j}{a}\right)\right).
\end{multline*}
We again use $a=2$.

\paragraph{Tests based on powers and minima.}
The test of \citet{allison2022distribution} compares empirical distribution functions for powers of the observations and minima of $m$ observations. Their difference is
\[ \Delta_{n,m} (x) := \frac{1}{n} \sum_{j=1}^{n} \mathbf{1}_{\{X_j^{1/m}\leq x\}} -\frac{1}{n^m}  \sum_{j_1,j_2,\ldots,{j_m}=1}^{n} \mathbf{1}_{\{\min(X_{j_1},X_{j_2},\ldots,{X_{j_m}})\leq x\}}.  \]
Integrating this difference gives the statistic
\[I_{n,m} :=\int_{1}^{\infty} \Delta_{n,m} (x)\,\mathrm{d}\widehat{F}_n(x) .
\]
We use the tuning parameter $m=2$.

\noindent The four statistics from \citet{ndwandwe2023new} use two sets of rank weights. Define
\[v_{j,m}:= \frac{1}{n^m} \bigg((n-j+1)^m -(n-j)^m\bigg) \text{ and }u_{j,m} := \binom{n-j}{m-1} . \]
The weights $v_{j,m}$ are used for $1\le j\le n$, and the weights $u_{j,m}$ for $1\le j\le n-m+1$. The first pair of statistics uses $v_{j,m}$:
\begin{multline*} S_{n,m,a}^{(1)} := \frac{1}{n} \sum_{j=1}^{n} \sum_{k=1}^{n} \left( \frac{2a}{a^2+( X_{(j)}^{1/m} -X_{(k)}^{1/m} )^2} -nv_{j,m} \frac{4a}{a^2+( X_{(j)} -X_{(k)}^{1/m} )^2}\right. \\  \left.
+ n^2 v_{j,m}v_{k,m} \frac{2a}{a^2+( X_{(j)} -X_{(k)})^2} \right), \end{multline*}
and
\begin{multline*}
	 S_{n,m,a}^{(2)} :=  \frac{1}{n} \sqrt{\frac{\pi}{a}} \sum_{j=1}^{n} \sum_{k=1}^{n}  \bigg( \exp \left( \frac{-(X_{(j)}^{1/m} -X_{(k)}^{1/m} )^2}{4a} \right)\\ -2nv_{j,m} \exp \left( \frac{-( X_{(j)} -X_{(k)}^{1/m} )^2}{4a} \right) \\
+ n^2 v_{j,m}v_{k,m} \exp \left( \frac{-( X_{(j)} -X_{(k)})^2}{4a} \right) \bigg).
\end{multline*}
The second pair uses $u_{j,m}$:
\begin{multline*}
 T_{n,m,a}^{(1)} := \frac{1}{n} \sum_{j=1}^{n} \sum_{k=1}^{n}  \frac{2a}{a^2+( X_{(j)}^{1/m} -X_{(k)}^{1/m} )^2} \\ -\sum_{j=1}^{n-m+1} \sum_{k=1}^{n} \binom{n}{m}^{-1} u_{j,m} \frac{4a}{a^2+( X_{(j)} -X_{(k)}^{1/m} )^2} \\
 + n \sum_{j=1}^{n-m+1} \sum_{k=1}^{n-m+1} \binom{n}{m}^{-2} u_{j,m} u_{k,m} \frac{2a}{a^2+( X_{(j)} -X_{(k)})^2},
\end{multline*}
and
\begin{multline*}
	T_{n,m,a}^{(2)} := \sqrt{\frac{\pi}{a}} \frac{1}{n}  \sum_{j=1}^{n} \sum_{k=1}^{n}  \exp \left( \frac{-(X_{(j)}^{1/m} -X_{(k)}^{1/m} )^2}{4a} \right) \\
 -2 \sqrt{\frac{\pi}{a}} \sum_{j=1}^{n-m+1} \sum_{k=1}^{n} \binom{n}{m}^{-1} u_{j,m} \exp \left( \frac{-( X_{(j)} -X_{(k)}^{1/m} )^2}{4a} \right) \\
	+ n \sqrt{\frac{\pi}{a}} \sum_{j=1}^{n-m+1} \sum_{k=1}^{n-m+1}  \binom{n}{m}^{-2} u_{j,m}u_{k,m} \exp \left( \frac{-( X_{(j)} -X_{(k)})^2}{4a} \right).
\end{multline*}
Within each pair, superscripts $(1)$ and $(2)$ distinguish the rational and Gaussian terms, respectively. We use $m=3$ and $a=2$ for all four statistics. Their asymptotic properties under the conditions imposed in the original work are discussed by \citet{ndwandwe2023new}.

\paragraph{Tests with fitted-shape weights.}
The two classes proposed by \citet{Ngatchou} incorporate the fitted shape parameter through the factors $X_j^{-\widehat\alpha_n(m-1)}$. To display their common structure, take $a>0$ and define the kernels
\[
 H_{1,a}(u):=\frac{2a}{a^2+u^2},\qquad
 H_{2,a}(u):=\sqrt{\frac{\pi}{a}}\exp\left(-\frac{u^2}{4a}\right),
 \qquad u\in\mathbb R.
\]
The two statistics can then be written as
\begin{align*}
 N^{(r)}_{n,m,a}
 &:=\frac1n\sum_{j,k=1}^n\bigg[
 \frac1{m^2}H_{r,a}\left(X_j^{1/m}-X_k^{1/m}\right)\\
 &\quad-\frac1m X_k^{-\widehat\alpha_n(m-1)}
 H_{r,a}\left(X_j^{1/m}-X_k\right)
 -\frac1m X_j^{-\widehat\alpha_n(m-1)}
 H_{r,a}\left(X_k^{1/m}-X_j\right)\\
 &\quad+X_j^{-\widehat\alpha_n(m-1)}X_k^{-\widehat\alpha_n(m-1)}
 H_{r,a}(X_j-X_k)\bigg],\qquad r=1,2.
\end{align*}
Following \citet{Ngatchou}, we use $(a,m)=(0.5,2)$ and $(0.5,3)$ for each class.

\noindent The test based on $I_{n,m}$ rejects for large absolute values, and the test based on $O_n$ uses both tails of its null distribution. The remaining competing tests reject for large values of their statistics.

To compare power, we draw samples from the alternatives in Table~\ref{AD}, which include distributions used in earlier simulation studies of Pareto goodness-of-fit tests. Distributions originally supported on $[0,\infty)$ are shifted by one so that all alternatives have support $[1,\infty)$. All densities in Table~\ref{AD} are zero for $x<1$, and all parameters used in the study are positive. The notation $\Gamma(\theta)$ denotes a gamma distribution with shape $\theta$ and rate one, shifted by one. The label $\mathrm{LogW}(\theta)$ is retained to match the simulation files: we generate $X:=\exp(Y)$, where $Y$ has a gamma distribution with shape and rate both equal to $1+\theta$. This gives the scaled log-gamma density in Table~\ref{AD}. The additional alternative $F_{1,3}$ is defined in Proposition~\ref{prop:continuous_invisible}.

Tables~\ref{power20}, \ref{power50} and \ref{power100} report rejection percentages for the three sample sizes. The first three rows of each table give empirical size under the Pareto model; the remaining rows give power against the alternatives. The final three columns contain the proposed procedures. The results support the following conclusions.
\begin{itemize}
\item[(a)] The simulated type-I error rates are generally near the nominal level, with visible finite-sample deviations. For example, under $\mathcal P(2)$, $S_{n,K}$ rejects in approximately $3.9\%$, $4.2\%$ and $4.8\%$ of samples for $n=20,50,100$, respectively. The fixed numerical calibration is therefore somewhat conservative in this case at smaller sample sizes.
\item[(b)] No procedure dominates across the alternatives considered.
\item[(c)] The proposed tests are competitive for several alternatives. However, $T_n$ is not omnibus, as Proposition~\ref{prop:continuous_invisible} shows. For $F_{1,3}$ its rejection percentages remain near $5\%$. Although $S_n$ is consistent under the stated assumptions, both calibrations also have low power against this alternative at the sample sizes studied.
\end{itemize}

\section{Direct comparison of bootstrap and covariance estimation methods}
\label{sec:Sn_comparison}
\begin{figure}[htbp]
\centering
\includegraphics[width=\linewidth]{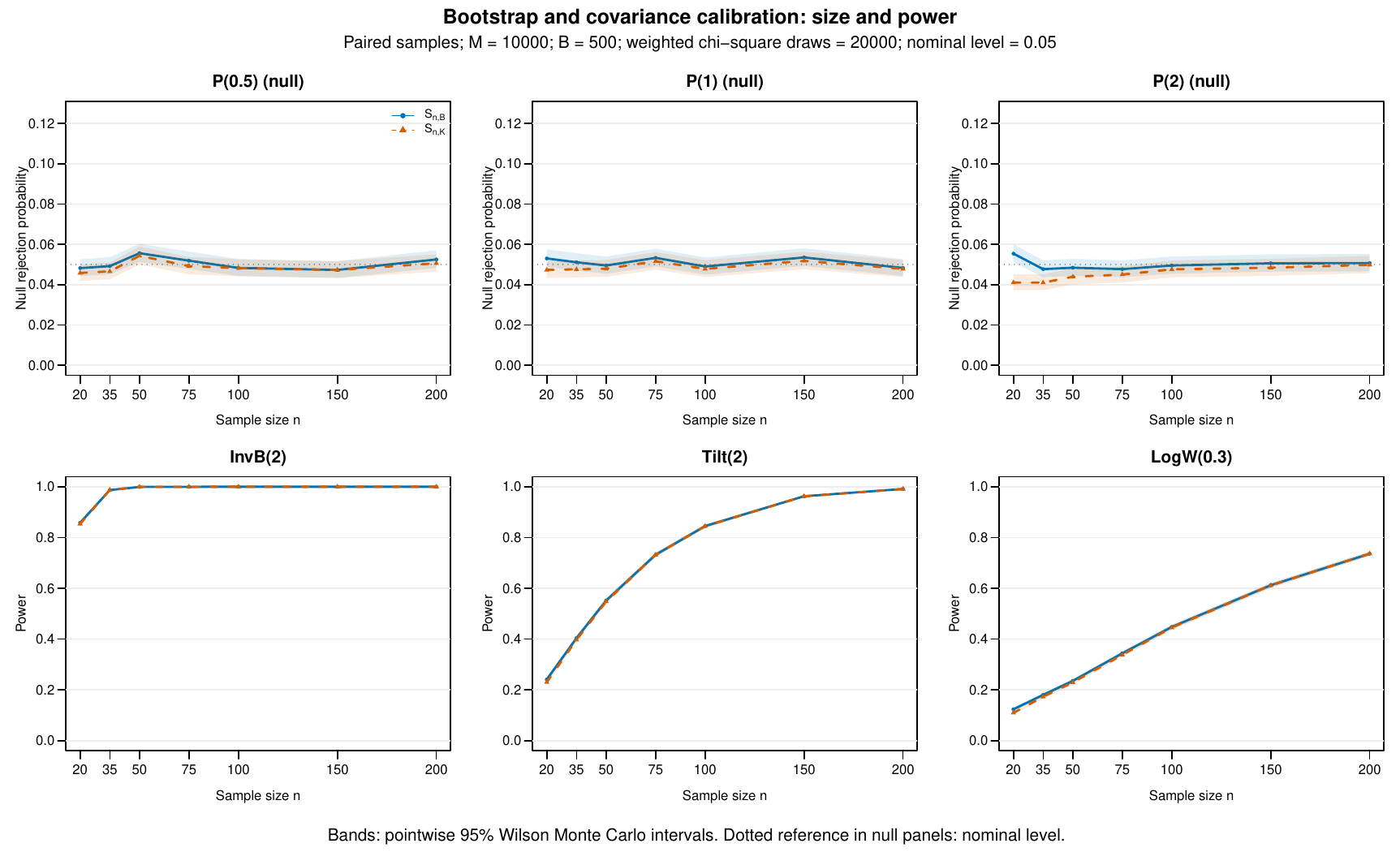}
\caption{Empirical size and power of $S_{n,B}$ (blue solid lines) and $S_{n,K}$ (orange dashed lines), based on $10{,}000$ paired samples per setting. Upper panels: Pareto null distributions. Lower panels: alternatives from Table~\ref{AD}; $\mathrm{LogW}(0.3)$ denotes the scaled log-gamma alternative. Shading gives pointwise $95\%$ Wilson Monte Carlo intervals for each rejection probability; dotted horizontal lines mark the nominal level $0.05$.}
\label{fig:Sn_power}
\end{figure}

\begin{figure}[htbp]
\centering
\includegraphics[width=\linewidth]{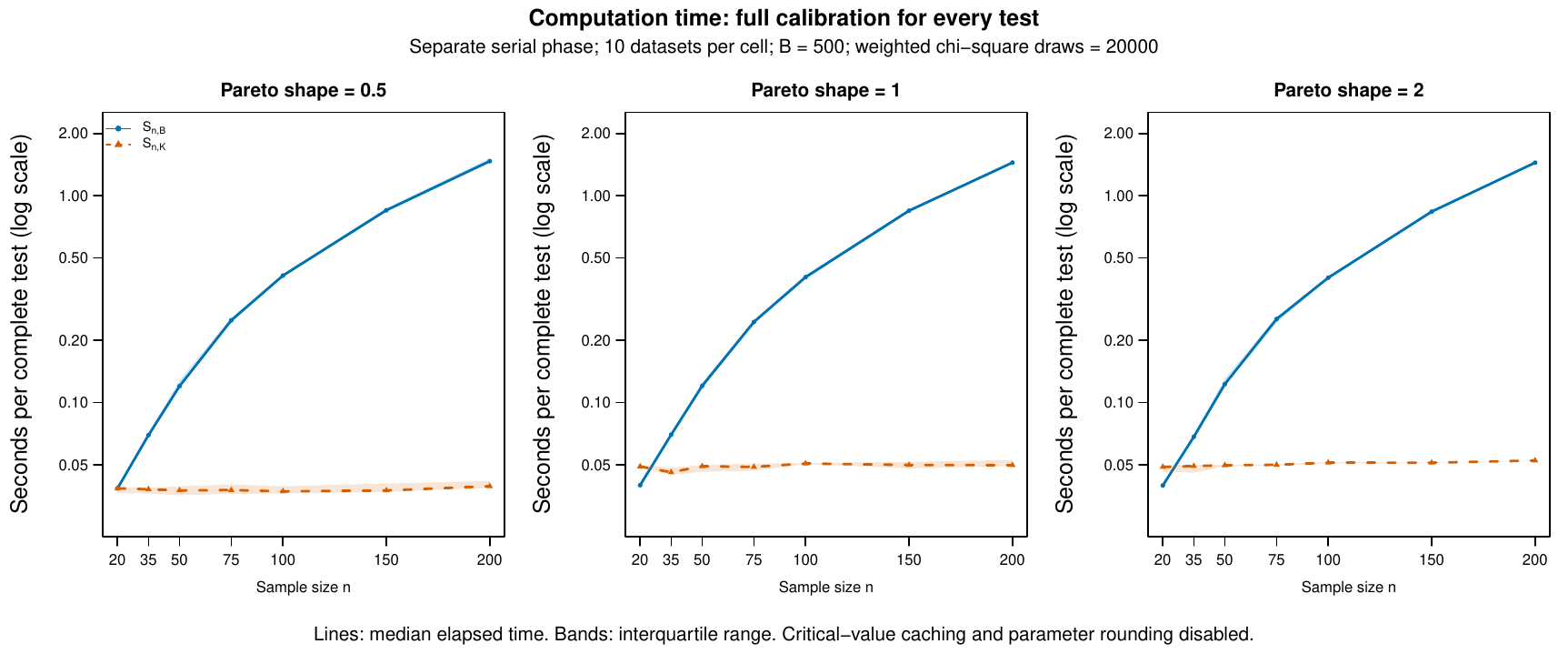}
\caption{Elapsed time per complete test under $\mathcal P(\alpha)$, measured serially on ten paired data sets per $(n,\alpha)$ setting. Lines show medians and shading shows interquartile ranges on a logarithmic time scale. Timings include fitting, the observed statistic, and full calibration, but exclude generation of the observed sample and warm-up calls. Fast calls are batched for at least $0.1$ seconds and averaged within each data set.}
\label{fig:Sn_time}
\end{figure}

For the distributions in Figure~\ref{fig:Sn_power}, we compare $S_{n,B}$ and $S_{n,K}$ using $M=10{,}000$ paired samples at each sample size $n=20$, $35$, $50$, $75$, $100$, $150$, and $200$. At nominal level $5\%$, we use $B=500$ bootstrap samples with parameter re-estimation, $20{,}000$ weighted chi-square draws, and the quadrature settings of Section~\ref{sec:Sn_covariance}. We used 14 parallel workers for the power experiment and
measured computation times separately in a single worker. Figure~\ref{fig:Sn_power} shows similar rejection rates, although $S_{n,K}$ is somewhat conservative under $\mathcal P(2)$ at the smaller sample sizes. The largest observed power difference is $1.42$ percentage points, in favor of $S_{n,B}$ under $\mathrm{LogW}(0.3)$ at $n=20$. The computational difference is more pronounced (see Figure~\ref{fig:Sn_time}): over the sample sizes studied, median covariance-calibrated test times remain near $0.04$--$0.05$ seconds, whereas bootstrap times rise from about $0.04$ to $1.45$ seconds. Covariance calibration is faster from $n=35$ onward and about $28$--$37$ times faster at $n=200$. Thus, for these settings and alternatives, it substantially reduces computation at larger sample sizes while giving similar power.

\section{Data analysis: Danish fire-insurance losses}
\label{datanalysis}

We consider the Danish fire-insurance losses available as \texttt{danish} in the R package \texttt{evir}
\citep{evir2018}. The data set contains $2{,}167$ claims recorded between 1980 and 1990. The amounts are expressed in millions of Danish kroner (DKK), adjusted to 1985 prices. These data are a standard benchmark in extreme-value analysis; see \citet{mcneil1997tails}.

Our aim is to examine how the assessment of an exact Pareto type-I model changes with the lower loss
threshold. Let $Y$ denote a recorded loss and consider the conditional hypothesis
$H_0(u):\mathbb P(Y>y\mid Y>u)=(y/u)^{-\alpha_u}$ for $y\geq u$ and some $\alpha_u>0$.
For each threshold $u$, we retain losses strictly exceeding $u$ and normalize them as
$X_j^{(u)}=Y_j/u$. The resulting testing problem has the known lower endpoint $1$, as assumed in
Section~\ref{Fusion}.

To assess threshold sensitivity, we use the grid $u\in\{1,1.5,2,3,5,10\}$, measured in millions of DKK, and retain only losses strictly exceeding each threshold. At $u=1$, excluding the 11 claims recorded exactly at the cutoff leaves $n_u=2{,}156$ exceedances. At $u=1.5$, six claims equal the cutoff and $1{,}386$ exceed it.

At each threshold, we estimate the Pareto shape by $\widehat\alpha_u=
\{n_u^{-1}\sum_{j=1}^{n_u}\log X_j^{(u)}\}^{-1}$. We calculate the signed Stein--Laplace statistic $T_{n_u}$ and the integrated squared statistic $n_uS_{n_u}$.

The signed statistic is calibrated by a two-sided parametric bootstrap. For $S_{n_u}$, we use both the bootstrap calibration and the plug-in covariance-operator calibration developed in Section~\ref{sec:Sn_covariance}. The bootstrap uses $B=4{,}999$ samples of size $n_u$ from $\mathcal P(\widehat\alpha_u)$, with shape re-estimation in every sample. The covariance-operator calculation uses the truncation interval $[0,12]$, $50$ Gauss--Legendre nodes, $80$ Gauss--Laguerre nodes, and $M=200{,}000$ weighted chi-square draws. Both calibrations use the same observed statistic $n_uS_{n_u}$.

For comparison, we also consider the Kolmogorov--Smirnov and Anderson--Darling statistics and the three Stein--Laplace procedures of \citet{bhati2026laplace}. In terms of their empirical discrepancy $R_{n_u}$, described in Section~\ref{Fusion}, these are $\mathscr{DS}_1=\int_0^\infty R_{n_u}(t)\,\mathrm dt$, $\mathscr{DS}_2=\sup_{t>0}|R_{n_u}(t)|$, and $\mathscr{DS}_3=\int_0^\infty R_{n_u}(t)^2\,\mathrm dt$. The implementation evaluates $\mathscr{DS}_1$ in closed form, approximates $\mathscr{DS}_2$ on a numerical grid, and evaluates $\mathscr{DS}_3$ by quadrature; the same rules are applied to the observed and bootstrap samples.

For a right-tailed statistic with observed value $r$ and bootstrap values $r_1^*,\ldots,r_B^*$, we report $\widehat p_B=(1+\sum_{b=1}^B\mathbf{1}_{\{r_b^*\geq r\}})/(B+1)$. The signed statistics $T_{n_u}$ and $\mathscr{DS}_1$ use twice the smaller plus-one tail probability, capped at one; the remaining statistics use the upper tail. Operator-calibrated $p$-values use the analogous correction with $M$ reference draws.

Table~\ref{tab:danish-thresholds} reports the numbers of exceedances, estimated shape parameters, test-statistic values, and calibrated $p$-values.

\begin{table}[t]
\centering
\caption{Threshold sensitivity of Pareto goodness-of-fit tests for the Danish fire-insurance losses. Thresholds are measured in millions of DKK. Bootstrap $p$-values use $B=4{,}999$ replications; the operator-calibrated $p$-values use $M=200{,}000$ weighted chi-square draws. All analyses use strict exceedances $Y>u$.}
\label{tab:danish-thresholds}
\small
\renewcommand{\arraystretch}{1.10}
\setlength{\tabcolsep}{5pt}
\begin{tabular*}{\linewidth}{@{\extracolsep{\fill}}lrrrrrr@{}}
\toprule
Threshold $u$ & 1.0 & 1.5 & 2.0 & 3.0 & 5.0 & 10.0 \\
\midrule
$n_u$ & 2156 & 1386 & 903 & 532 & 254 & 109 \\
$\widehat\alpha_u$ & 1.264 & 1.406 & 1.371 & 1.421 & 1.414 & 1.614 \\
$T_{n_u}$ & 0.08269 & -0.00195 & 0.02342 & 0.01134 & 0.00553 & 0.05758 \\
$n_uS_{n_u}$ & 4.40989 & 0.00198 & 0.16338 & 0.02069 & 0.02676 & 0.11037 \\
\midrule
$p_{T,B}$ & 0.0004 & 0.8892 & 0.2516 & 0.7040 & 0.9340 & 0.3880 \\
$p_{S,B}$ & 0.0002 & 0.9462 & 0.2326 & 0.6908 & 0.6340 & 0.3442 \\
$p_{S,K}$ & $5.0\cdot10^{-6}$ & 0.9475 & 0.2363 & 0.6916 & 0.6451 & 0.3223 \\
$p_{KS}$ & 0.0002 & 0.5920 & 0.3444 & 0.4500 & 0.1984 & 0.5684 \\
$p_{AD}$ & 0.0002 & 0.3634 & 0.3996 & 0.4678 & 0.1340 & 0.5174 \\
$p_{\mathscr{DS}_1}$ & 0.0004 & 0.8548 & 0.2608 & 0.7548 & 0.8620 & 0.3816 \\
$p_{\mathscr{DS}_2}$ & 0.0002 & 0.9576 & 0.2080 & 0.5850 & 0.4934 & 0.4544 \\
$p_{\mathscr{DS}_3}$ & 0.0002 & 0.9238 & 0.2474 & 0.6952 & 0.5558 & 0.3570 \\
\bottomrule
\end{tabular*}

\par\smallskip
\begin{minipage}{0.94\linewidth}
\footnotesize
At $u=1$, the probabilities $0.0002$ and $5.0\cdot10^{-6}$ are the respective plus-one Monte Carlo reporting minima for $B=4{,}999$ and $M=200{,}000$; no reference value exceeded the observed statistic. The minimum for the doubled-tail rule is $0.0004$.
\end{minipage}
\end{table}

At the reporting threshold $u=1$, the fitted shape parameter is $\widehat\alpha_u=1.264$, and all procedures reject the
Pareto model at the $5\%$ level. At $u=1.5$, neither calibration of $S_{n_u}$ rejects, although
$1{,}386$ exceedances remain. The bootstrap and operator $p$-values are $0.9462$ and $0.9475$, respectively.
The Kolmogorov--Smirnov and Anderson--Darling tests also do not reject, with $p$-values $0.5920$ and $0.3634$.
The three related Stein--Laplace procedures reach the same
conclusion.

At the remaining thresholds, neither calibration of $S_{n_u}$ rejects at the nominal $5\%$ level. The fitted shape parameter varies between approximately $1.37$ and $1.42$ for thresholds from $1.5$ to $5$ million DKK. At $u=10$,
the estimate increases to $1.61$, although the number of exceedances has declined to $109$.

Across the six thresholds, bootstrap and covariance-operator calibration yield the same rejection decisions and generally
similar $p$-values. Both detect a departure from the Pareto model at the reporting threshold of one million DKK, whereas
neither rejects at the five higher thresholds considered. The substantive conclusions of this application are therefore
unchanged by the choice between the two calibration methods.

\section{Conclusion} \label{conclusions}
The central distinction between the proposed tests is the information they retain from the Stein--Laplace characterization. The signed statistic $T_n$ reduces the discrepancy to a scalar moment restriction, yielding a simple diagnostic with an explicit studentization and a standard normal limiting null distribution. However, continuous non-Pareto distributions can satisfy the same restriction, so the resulting test is not consistent against every alternative. In contrast, squaring before integration preserves the characterizing property. The statistic $S_n$ converges to a population discrepancy that is strictly positive for every non-Pareto distribution on $[1,\infty)$ satisfying $0<\mathbb{E}_F[\log X]<\infty$. Together, these results show why the manner in which a characterization is aggregated is decisive for the consistency of the resulting test.

The covariance structure of the limiting Gaussian process provides both a description of the null law of $nS_n$ and a practical alternative to bootstrap calibration. Approximating the eigenvalues of the fitted covariance operator allows critical values to be obtained from a weighted chi-square distribution without repeatedly generating and refitting Pareto samples. The direct comparison of the two calibrations demonstrates the computational value of this approach. In the reported timing experiment, covariance-operator calibration was approximately $28$--$37$ times faster at $n=200$, while the paired power experiment showed similar rejection probabilities for the alternatives examined. The parametric bootstrap nevertheless remains useful in small samples, where the operator approximation was conservative in some settings. Its occasionally higher empirical power must therefore be considered together with the corresponding differences in empirical size.

The broader comparison with existing procedures does not support a universal ranking. The proposed tests are competitive against several alternatives, but established procedures are appreciably more powerful in other cases, a pattern consistent with the theoretical findings of \citet{janssen2000global}. The constructed non-Pareto mixture is particularly informative: it demonstrates the lack of consistency of $T_n$, while also exposing low finite-sample power of both calibrations of $S_n$.

The Danish fire-insurance application demonstrates the importance of specifying the region over which an exact Pareto model is
assessed. The model is strongly rejected at the original reporting threshold of one million DKK, whereas neither calibration of
$S_n$ rejects at any of the five higher thresholds considered. In particular, the nonrejection at $1.5$ million DKK occurs
with $1{,}386$ exceedances remaining. These findings are compatible with a Pareto approximation farther into the
upper tail, but they neither identify a statistically validated threshold nor establish exact Pareto behavior. Developing
selection-aware inference for the Stein--Laplace statistics therefore remains a natural extension of the present theory.

Several further questions arise directly from the results. The known-endpoint assumption suggests two distinct extensions:
estimation of the support endpoint within an exact Pareto model, and selection of a threshold above which a Pareto approximation is adequate. Goodness-of-fit testing with an estimated endpoint has been studied by \citet{ndwandwe2023testing}, while procedures combining threshold selection and Pareto fitting include those of
\citet{grama2008statistics} and \citet{clauset2009power}. The asymptotic analysis of \citet{drees2020minimum} shows that
threshold selection can alter the limiting distribution of the subsequent tail-index estimator. For the present Stein--Laplace tests, the corresponding task is to establish the joint limit of the selected threshold, the shape estimator, and the empirical discrepancy, and to determine when covariance-operator or bootstrap calibration remains valid. Ordered testing procedures \citep{bader2018automated} and methods that propagate threshold-selection uncertainty \citep{murphy2025automated} provide relevant starting points for this extension. Another interesting line of further research includes limits under local alternatives and asymptotic efficiency calculations. They would clarify which departures are weakly represented by the unweighted Stein--Laplace discrepancy. Such an analysis could guide the choice of weighted versions of $S_n$, with the aim of improving sensitivity while retaining consistency.

\section*{Declaration of generative AI and AI-assisted technologies in the manuscript preparation process}

OpenAI's ChatGPT 6 Astra was used as an AI-assisted research tool in the
preparation of this work. In particular, the code for the simulation study was partially written by ChatGPT. All suggested proofs were checked and independently verified by the authors. ChatGPT was also used for assistance with exposition, \LaTeX{} formulation, and consistency checks. The authors take full responsibility for all results and conclusions in the paper.

\bibliographystyle{apalike}
\bibliography{bibitex}

\begin{appendix}
\section{Tables}
We provide the tables of alternatives and simulation results of Section \ref{MCstudy}.
\begin{table}[htbp]
\begin{center}
\caption{Alternative distributions and their densities on $x>1$.}
\label{AD}
\begin{tabular}{ l  c c }
\hline
Alternative & Probability density function & Notation \\
\hline
Gamma& $\frac{1}{\Gamma(\theta)} (x-1)^{\theta-1} \exp\left(-(x-1)\right)$ & $\Gamma(\theta)$ \\
Inverse-beta & $\frac{1+\theta}{x^2} \left(1- \frac{1}{x}\right)^ \theta$ & $\mathrm{InvB}(\theta)$ \\
Tilted Pareto & $\frac{1+\theta}{(x+\theta)^2}$ & $\mathrm{Tilt}(\theta)$ \\
Benini & $\frac{\exp\left(-\theta \log^2 x\right)}{x^2} (1+2\theta \log x)$ & $\mathrm{Ben}(\theta)$ \\
Log-gamma & $\frac{\log^\theta x}{x^2 \Gamma (1+\theta)}$ & $\mathrm{LogG}(\theta)$\\
Log-normal& $\frac{\exp \left(-\frac{1}{2} \left(\frac{\log (x-1)}{\theta}\right)^2\right)}{\theta (x-1) \sqrt{2 \pi}}$ & $\mathrm{LN}(\theta)$ \\
Rayleigh & $\frac{x-1}{\theta^2} \exp\left(-\frac{(x-1)^2}{2\theta^2}\right)$ & $\mathrm{Ray}(\theta)$\\
Weibull & $\theta (x-1)^{\theta-1} \exp\left(-(x-1)^\theta\right)$ & $\mathrm{W}(\theta)$\\
L\'evy & $\sqrt{\frac{\theta}{2\pi}} \cdot \frac{\exp\left(-\frac{\theta}{2(x-1)}\right)}{(x-1)^{3/2}}$ & $\mathrm{L}(\theta)$ \\
Burr & $cb \frac{\left(\frac{x-1}{a}\right)^{b-1}}{a \left(1+ {(\frac{x-1}{a})}^b\right) ^{c+1}}$ &$\mathrm{Burr}(a,b,c)$ \\
Inverse-Gaussian & $\sqrt{\frac{\lambda}{2\pi (x-1)^3}} \exp\left(-\frac{\lambda(x-1-\mu)^2}{2\mu^2(x-1)}\right)$ & $\mathrm{IG}(\mu,\lambda)$\\
Scaled log-gamma & $\frac{(1+\theta)^{1+\theta}\log^\theta x}{\Gamma(1+\theta)x^{2+\theta}}$ & $\mathrm{LogW}(\theta)$\\
Fr\'echet & $\frac{a}{b} \left(\frac{x-1}{b}\right)^{-(a+1)} \exp \left( - \left(\frac{x-1}{b}\right)^{-a}\right)$ & $\mathrm{FR}(a,b)$\\
Half-normal & $\sqrt{\frac{2}{\pi \theta^2}} \exp\left(-\frac{(x-1)^2}{2\theta^2}\right)$ & $\mathrm{HN}(\theta)$\\
Chi-square & $\frac{1}{2^{\nu/2} \Gamma (\nu/2)} (x-1)^{\nu/2-1} \exp\left(-(x-1)/2\right)$ & $\chi^2_\nu$ \\
Dhillon & $\frac{\theta+1}{x} \exp \left( -(\log(x) )^{\theta+1} \right) (\log(x))^ \theta$ & $\mathrm{DH}(\theta)$ \\
Log-logistic & $\frac{\gamma \left(\frac{x-1}{\theta}\right)^\gamma}{(x-1)[1+\left(\frac{x-1}{\theta}\right)^\gamma]^2}$ & $\mathrm{LLogis}(\gamma,\theta)$\\
Linear failure rate & $\left(1+\theta(x-1)\right)\exp\left(-(x-1)-\theta(x-1)^2/2\right)$ & $\mathrm{LFR}(\theta)$\\ \hline
\end{tabular}
\end{center}
\end{table}

\begin{sidewaystable}[t]
\caption{Empirical size and power (percentages, rounded to one decimal place) of the proposed and competing tests ($n=20$, $N=10{,}000$, $B=500$, $\eta=0.05$). The distribution $F_{1,3}$ is the alternative in Proposition~\ref{prop:continuous_invisible} with $\lambda=1$ and $b=3$.} \label{power20}
\resizebox{\linewidth}{!}{
\begin{tabular}{|l|llllllllllllllllll|} \hline
Distribution  & $O_n$ & $KS_n$ & $CM_n$ & $AD_n$ & $G_{n,2}$ & $M_{n,2}$ & $|I_{n,2}|$ & $S_{n,3,2}^{(1)}$ & $S_{n,3,2}^{(2)}$ & $T_{n,3,2}^{(1)}$ & $T_{n,3,2}^{(2)}$ & $N^{(1)}_{n,2,0.5}$ & $N^{(1)}_{n,3,0.5}$ & $N^{(2)}_{n,2,0.5}$ & $N^{(2)}_{n,3,0.5}$ & ${T}_n$ & $S_{n,K}$ & $S_{n,B}$ \\ \hline
$\mathcal P(0.5)$ & 5.1 & 4.9 & 5.2 & 5.2 & 5.1 & 5.2 & 4.8 & 4.4 & 4.3 & 4.2 & 4.0 & 5.4 & 5.2 & 5.2 & 5.0 & 5.3 & 4.8 & 5.0 \\
$\mathcal P(1)$ & 5.1 & 5.2 & 5.1 & 5.0 & 5.3 & 5.0 & 4.8 & 4.0 & 3.4 & 3.3 & 3.0 & 5.1 & 4.8 & 4.6 & 5.0 & 5.4 & 4.7 & 5.2 \\
$\mathcal P(2)$ & 5.1 & 5.0 & 5.1 & 5.2 & 5.4 & 4.9 & 4.9 & 3.3 & 3.2 & 3.3 & 3.2 & 4.6 & 4.8 & 4.7 & 4.7 & 5.5 & 3.9 & 5.4 \\
\hline
$\Gamma(0.5)$ & 32.2 & 21.3 & 23.8 & 45.2 & 13.7 & 19.6 & 30.0 & 6.0 & 5.0 & 17.2 & 15.7 & 20.1 & 27.1 & 16.4 & 25.5 & 13.0 & 10.1 & 12.5 \\
$\Gamma(1.2)$ & 54.6 & 45.2 & 55.9 & 50.3 & 60.5 & 60.6 & 50.4 & 62.1 & 61.7 & 54.4 & 52.3 & 46.7 & 53.1 & 59.9 & 58.1 & 59.8 & 57.6 & 62.2 \\
$\Gamma(2)$ & 99.4 & 96.0 & 99.3 & 99.1 & 99.7 & 99.7 & 98.7 & 99.5 & 99.6 & 99.1 & 99.2 & 95.4 & 98.5 & 99.2 & 99.6 & 99.6 & 99.7 & 99.7 \\
$\mathrm{InvB}(0.1)$ & 5.3 & 6.0 & 5.8 & 4.9 & 6.1 & 5.8 & 5.3 & 5.7 & 5.1 & 4.8 & 4.0 & 6.1 & 5.9 & 5.6 & 5.9 & 6.0 & 5.4 & 6.2 \\
$\mathrm{InvB}(2)$ & 75.0 & 69.1 & 77.3 & 76.0 & 84.5 & 77.7 & 77.0 & 72.8 & 72.3 & 70.3 & 69.9 & 41.8 & 63.6 & 54.5 & 76.3 & 84.9 & 84.7 & 85.2 \\
$\mathrm{Tilt}(1)$ & 11.5 & 11.3 & 12.5 & 10.3 & 12.7 & 12.8 & 11.8 & 14.0 & 13.6 & 11.9 & 11.4 & 9.0 & 9.4 & 11.1 & 12.4 & 11.4 & 11.9 & 12.7 \\
$\mathrm{Tilt}(2)$ & 23.9 & 22.1 & 25.8 & 21.3 & 25.7 & 26.6 & 23.6 & 25.8 & 26.5 & 23.5 & 23.3 & 13.2 & 16.1 & 18.2 & 23.6 & 22.3 & 24.1 & 25.0 \\
$\mathrm{Tilt}(3)$ & 35.6 & 31.5 & 38.1 & 32.3 & 36.3 & 39.3 & 34.2 & 34.7 & 36.5 & 31.7 & 32.9 & 15.1 & 20.4 & 23.7 & 32.0 & 30.6 & 33.5 & 34.5 \\
$\mathrm{Ben}(0.1)$ & 6.1 & 6.2 & 6.5 & 5.6 & 6.3 & 6.1 & 6.1 & 6.7 & 6.1 & 5.4 & 4.4 & 6.7 & 6.1 & 6.5 & 6.1 & 6.1 & 5.4 & 6.7 \\
$\mathrm{Ben}(0.3)$ & 10.7 & 11.2 & 12.0 & 9.7 & 12.0 & 12.1 & 10.5 & 13.6 & 12.9 & 10.8 & 9.4 & 10.4 & 11.1 & 12.1 & 11.9 & 11.8 & 10.7 & 12.8 \\
$\mathrm{Ben}(0.7)$ & 20.6 & 19.2 & 22.7 & 18.8 & 23.5 & 23.8 & 19.6 & 26.0 & 24.9 & 20.7 & 18.9 & 19.4 & 21.7 & 22.7 & 22.8 & 23.8 & 20.6 & 24.6 \\
$\mathrm{LogG}(0.5)$ & 19.3 & 17.3 & 20.3 & 17.3 & 22.7 & 20.8 & 19.4 & 20.2 & 19.8 & 18.1 & 17.3 & 12.7 & 15.8 & 15.2 & 20.6 & 21.9 & 22.0 & 23.0 \\
$\mathrm{LogG}(1.2)$ & 57.9 & 50.1 & 59.7 & 56.9 & 67.6 & 62.3 & 58.4 & 52.2 & 52.5 & 49.8 & 49.9 & 23.7 & 37.1 & 34.5 & 53.2 & 64.9 & 65.5 & 66.5 \\
$\mathrm{LogG}(2)$ & 88.3 & 81.3 & 89.1 & 88.4 & 93.5 & 90.6 & 88.4 & 77.8 & 78.4 & 75.7 & 76.2 & 26.7 & 48.9 & 44.8 & 72.7 & 90.5 & 90.8 & 91.2 \\
$\mathrm{LN}(1)$ & 62.9 & 54.8 & 65.0 & 63.1 & 71.4 & 68.7 & 63.4 & 66.2 & 65.8 & 61.3 & 59.2 & 44.8 & 65.4 & 54.9 & 71.7 & 72.2 & 67.6 & 71.2 \\
$\mathrm{LN}(1.2)$ & 28.7 & 24.9 & 29.6 & 27.3 & 34.6 & 32.6 & 28.6 & 31.2 & 30.7 & 26.9 & 25.2 & 19.5 & 28.8 & 24.1 & 34.1 & 35.2 & 31.5 & 35.2 \\
$\mathrm{LN}(1.5)$ & 6.9 & 6.5 & 7.3 & 5.7 & 8.1 & 7.6 & 6.6 & 7.8 & 7.5 & 6.3 & 5.7 & 6.5 & 6.5 & 7.3 & 8.0 & 8.0 & 7.2 & 8.4 \\
$\mathrm{Ray}(1)$ & 100.0 & 99.5 & 100.0 & 100.0 & 100.0 & 100.0 & 99.9 & 100.0 & 100.0 & 100.0 & 100.0 & 99.9 & 100.0 & 100.0 & 100.0 & 100.0 & 100.0 & 100.0 \\
$\mathrm{W}(0.5)$ & 45.7 & 34.6 & 37.6 & 61.7 & 30.6 & 31.4 & 43.9 & 14.7 & 11.6 & 29.0 & 24.4 & 29.0 & 45.7 & 30.9 & 41.1 & 34.7 & 27.4 & 29.9 \\
$\mathrm{W}(1.2)$ & 59.5 & 49.2 & 60.8 & 55.5 & 65.7 & 65.4 & 55.2 & 67.3 & 66.4 & 59.8 & 57.4 & 54.9 & 60.0 & 65.8 & 62.8 & 65.8 & 61.3 & 67.1 \\
$\mathrm{W}(1.75)$ & 98.8 & 94.7 & 98.7 & 98.4 & 99.3 & 99.3 & 97.9 & 99.4 & 99.4 & 98.8 & 98.6 & 98.3 & 98.8 & 99.4 & 99.1 & 99.4 & 99.1 & 99.5 \\
$\mathrm{L}(2)$ & 21.1 & 20.4 & 22.7 & 22.1 & 30.1 & 20.9 & 24.8 & 26.1 & 23.3 & 25.5 & 22.2 & 18.1 & 25.6 & 17.7 & 27.4 & 34.8 & 31.0 & 32.4 \\
$\mathrm{Burr}(1.5,0.5,0.5)$ & 26.2 & 19.1 & 19.6 & 45.0 & 37.7 & 15.3 & 25.6 & 3.4 & 3.9 & 4.8 & 5.5 & 5.4 & 26.0 & 8.2 & 25.9 & 44.5 & 41.8 & 42.6 \\
$\mathrm{IG}(1,1)$ & 76.4 & 72.7 & 80.0 & 80.5 & 80.9 & 81.2 & 80.0 & 81.6 & 80.7 & 77.5 & 75.1 & 61.5 & 83.1 & 67.7 & 85.5 & 82.1 & 75.3 & 79.8 \\
$\mathrm{LogW}(0.1)$ & 5.5 & 5.5 & 5.9 & 4.9 & 6.2 & 5.9 & 5.5 & 5.5 & 5.2 & 4.7 & 4.2 & 5.8 & 5.8 & 5.4 & 6.1 & 6.1 & 5.6 & 6.3 \\
$\mathrm{LogW}(0.3)$ & 10.5 & 10.1 & 10.9 & 9.1 & 12.0 & 11.2 & 10.4 & 11.9 & 11.2 & 10.0 & 8.8 & 9.0 & 10.5 & 9.7 & 11.7 & 11.9 & 11.1 & 12.5 \\
$\mathrm{LogW}(1.75)$ & 81.9 & 73.7 & 82.8 & 81.4 & 87.4 & 85.2 & 81.9 & 82.2 & 81.9 & 79.4 & 78.2 & 60.9 & 80.8 & 71.9 & 86.5 & 87.8 & 86.3 & 87.7 \\
$\mathrm{FR}(1,1)$ & 47.1 & 48.3 & 51.1 & 51.6 & 51.9 & 44.0 & 55.4 & 48.6 & 42.4 & 47.0 & 39.4 & 37.8 & 59.3 & 34.4 & 53.0 & 56.4 & 51.2 & 53.2 \\
$\mathrm{HN}(0.8)$ & 57.8 & 48.7 & 59.1 & 53.5 & 61.6 & 62.0 & 52.1 & 65.2 & 63.5 & 56.8 & 54.8 & 62.5 & 57.0 & 63.0 & 56.3 & 63.5 & 55.3 & 63.0 \\
$\mathrm{HN}(1)$ & 63.2 & 53.5 & 65.0 & 59.5 & 67.0 & 66.8 & 57.6 & 70.4 & 68.9 & 63.0 & 60.3 & 65.8 & 61.9 & 70.6 & 61.9 & 67.7 & 63.4 & 68.7 \\
$\chi^2_4$ & 100.0 & 99.5 & 100.0 & 100.0 & 100.0 & 100.0 & 99.9 & 99.9 & 100.0 & 99.9 & 99.9 & 97.7 & 99.6 & 99.8 & 100.0 & 99.9 & 100.0 & 100.0 \\
$\mathrm{DH}(0.2)$ & 12.5 & 11.5 & 13.5 & 10.9 & 14.8 & 14.0 & 12.7 & 14.4 & 13.8 & 11.8 & 10.9 & 10.5 & 12.6 & 12.3 & 14.4 & 14.7 & 13.7 & 15.2 \\
$\mathrm{DH}(0.4)$ & 33.2 & 28.7 & 35.1 & 30.8 & 38.7 & 37.4 & 32.9 & 36.6 & 36.0 & 32.0 & 30.6 & 23.1 & 31.9 & 29.9 & 37.9 & 37.8 & 36.3 & 39.3 \\
$\mathrm{DH}(0.7)$ & 72.0 & 61.7 & 73.2 & 69.8 & 78.4 & 77.4 & 70.5 & 75.1 & 75.2 & 70.3 & 69.5 & 52.5 & 69.7 & 66.9 & 77.4 & 78.0 & 76.8 & 79.1 \\
$\mathrm{LLogis}(2,1)$ & 86.0 & 80.5 & 87.1 & 85.6 & 87.5 & 86.6 & 87.0 & 88.1 & 86.2 & 85.7 & 82.8 & 75.5 & 88.4 & 78.6 & 89.4 & 88.7 & 85.7 & 87.8 \\
$\mathrm{LFR}(0.2)$ & 36.9 & 31.3 & 38.8 & 33.7 & 41.4 & 41.5 & 33.0 & 44.4 & 43.0 & 36.8 & 34.0 & 37.5 & 36.5 & 43.8 & 38.0 & 41.9 & 37.2 & 43.0 \\
$\mathrm{LFR}(0.5)$ & 45.7 & 37.8 & 47.8 & 41.1 & 50.3 & 50.1 & 41.0 & 54.0 & 52.0 & 45.5 & 42.8 & 47.6 & 45.5 & 52.5 & 45.7 & 51.0 & 45.0 & 51.8 \\
$\mathrm{LFR}(0.8)$ & 50.4 & 42.3 & 52.4 & 46.2 & 54.6 & 55.1 & 45.7 & 58.4 & 56.8 & 50.0 & 47.6 & 53.8 & 50.1 & 56.3 & 50.0 & 56.1 & 48.7 & 56.2 \\
$F_{1,3}$ & 27.7 & 57.4 & 46.1 & 41.8 & 7.2 & 20.6 & 18.5 & 45.0 & 42.7 & 47.6 & 45.3 & 23.7 & 26.5 & 44.8 & 32.3 & 4.4 & 4.1 & 4.4 \\
\hline
\end{tabular}
}
\end{sidewaystable}

\begin{sidewaystable}[ht]
\caption{Empirical size and power (percentages, rounded to one decimal place) of the proposed and competing tests ($n=50$, $N=10{,}000$, $B=500$, $\eta=0.05$). The distribution $F_{1,3}$ is the alternative in Proposition~\ref{prop:continuous_invisible} with $\lambda=1$ and $b=3$.} \label{power50}
\resizebox{\linewidth}{!}{
\begin{tabular}{|l|llllllllllllllllll|} \hline
Distribution  & $O_n$ & $KS_n$ & $CM_n$ & $AD_n$ & $G_{n,2}$ & $M_{n,2}$ & $|I_{n,2}|$ & $S_{n,3,2}^{(1)}$ & $S_{n,3,2}^{(2)}$ & $T_{n,3,2}^{(1)}$ & $T_{n,3,2}^{(2)}$ & $N^{(1)}_{n,2,0.5}$ & $N^{(1)}_{n,3,0.5}$ & $N^{(2)}_{n,2,0.5}$ & $N^{(2)}_{n,3,0.5}$ & ${T}_n$ & $S_{n,K}$ & $S_{n,B}$ \\ \hline
$\mathcal P(0.5)$ & 5.3 & 4.9 & 5.0 & 5.1 & 5.4 & 5.2 & 4.9 & 4.7 & 4.7 & 4.7 & 4.5 & 4.7 & 4.9 & 5.0 & 5.2 & 5.2 & 5.1 & 5.3 \\
$\mathcal P(1)$ & 5.3 & 5.2 & 5.5 & 5.1 & 5.3 & 5.1 & 5.3 & 4.4 & 4.2 & 4.4 & 4.0 & 5.3 & 5.4 & 5.0 & 5.2 & 5.4 & 5.0 & 5.3 \\
$\mathcal P(2)$ & 5.0 & 5.0 & 4.9 & 4.7 & 4.7 & 4.7 & 5.1 & 4.1 & 3.9 & 4.0 & 3.9 & 4.8 & 4.9 & 4.9 & 4.9 & 4.9 & 4.2 & 4.8 \\
\hline
$\Gamma(0.5)$ & 62.0 & 49.1 & 54.0 & 78.1 & 24.0 & 40.3 & 63.4 & 32.5 & 26.9 & 44.8 & 39.6 & 44.7 & 56.4 & 32.2 & 52.8 & 23.7 & 19.1 & 21.0 \\
$\Gamma(1.2)$ & 93.6 & 86.3 & 94.3 & 93.7 & 96.9 & 96.7 & 91.7 & 96.7 & 97.3 & 95.2 & 95.9 & 89.4 & 93.0 & 96.9 & 95.3 & 96.3 & 96.9 & 97.2 \\
$\Gamma(2)$ & 100.0 & 100.0 & 100.0 & 100.0 & 100.0 & 100.0 & 100.0 & 100.0 & 100.0 & 100.0 & 100.0 & 100.0 & 100.0 & 100.0 & 100.0 & 100.0 & 100.0 & 100.0 \\
$\mathrm{InvB}(0.1)$ & 6.8 & 6.7 & 6.9 & 6.1 & 7.1 & 7.0 & 6.9 & 7.4 & 7.1 & 6.3 & 5.9 & 6.1 & 6.8 & 6.6 & 7.3 & 7.2 & 7.1 & 7.4 \\
$\mathrm{InvB}(2)$ & 99.7 & 99.0 & 99.7 & 99.8 & 99.9 & 99.6 & 99.8 & 98.7 & 98.6 & 98.4 & 98.2 & 86.5 & 98.7 & 94.3 & 99.5 & 99.9 & 99.9 & 99.9 \\
$\mathrm{Tilt}(1)$ & 24.4 & 22.1 & 25.9 & 22.7 & 27.0 & 27.9 & 23.7 & 24.4 & 25.9 & 22.1 & 22.8 & 13.0 & 17.4 & 18.5 & 24.3 & 23.9 & 26.1 & 26.4 \\
$\mathrm{Tilt}(2)$ & 53.9 & 47.4 & 56.3 & 51.9 & 56.3 & 59.6 & 52.4 & 48.5 & 51.4 & 45.3 & 48.2 & 23.0 & 35.5 & 35.6 & 50.0 & 50.4 & 54.2 & 54.7 \\
$\mathrm{Tilt}(3)$ & 75.5 & 68.3 & 77.3 & 73.7 & 76.0 & 80.3 & 73.4 & 66.2 & 69.7 & 63.3 & 66.6 & 31.6 & 51.1 & 49.5 & 68.6 & 68.3 & 73.1 & 73.1 \\
$\mathrm{Ben}(0.1)$ & 8.4 & 8.7 & 9.2 & 8.0 & 9.9 & 9.9 & 8.2 & 10.4 & 10.5 & 8.6 & 8.4 & 7.7 & 8.0 & 9.0 & 9.7 & 9.6 & 9.6 & 10.2 \\
$\mathrm{Ben}(0.3)$ & 23.0 & 21.3 & 25.9 & 22.5 & 28.7 & 29.0 & 21.9 & 29.7 & 31.5 & 25.1 & 25.8 & 18.6 & 22.6 & 26.7 & 26.9 & 27.7 & 28.1 & 29.4 \\
$\mathrm{Ben}(0.7)$ & 48.2 & 42.5 & 51.8 & 47.8 & 57.5 & 57.4 & 46.0 & 59.2 & 60.6 & 53.2 & 53.6 & 43.8 & 49.6 & 56.9 & 53.7 & 56.9 & 56.9 & 58.4 \\
$\mathrm{LogG}(0.5)$ & 45.4 & 36.7 & 44.0 & 44.8 & 51.4 & 46.3 & 46.8 & 39.0 & 39.6 & 36.2 & 36.5 & 20.9 & 36.2 & 28.6 & 44.8 & 52.2 & 51.4 & 51.8 \\
$\mathrm{LogG}(1.2)$ & 96.5 & 92.4 & 96.2 & 97.2 & 98.5 & 97.0 & 97.0 & 88.5 & 88.9 & 87.1 & 87.7 & 54.4 & 84.7 & 72.0 & 93.4 & 98.3 & 98.5 & 98.4 \\
$\mathrm{LogG}(2)$ & 100.0 & 99.8 & 100.0 & 100.0 & 100.0 & 100.0 & 100.0 & 99.2 & 99.1 & 99.0 & 99.0 & 73.1 & 97.4 & 89.9 & 99.5 & 100.0 & 100.0 & 100.0 \\
$\mathrm{LN}(1)$ & 98.6 & 95.4 & 98.5 & 99.2 & 99.3 & 99.0 & 98.9 & 98.3 & 98.5 & 97.7 & 97.9 & 90.3 & 98.5 & 95.1 & 99.2 & 99.2 & 99.1 & 99.1 \\
$\mathrm{LN}(1.2)$ & 69.9 & 57.4 & 68.4 & 72.5 & 76.8 & 74.2 & 71.5 & 68.0 & 71.0 & 63.9 & 65.8 & 43.6 & 67.5 & 57.1 & 75.9 & 77.1 & 75.5 & 76.0 \\
$\mathrm{LN}(1.5)$ & 11.9 & 9.8 & 12.0 & 11.8 & 15.4 & 14.6 & 11.8 & 12.6 & 13.8 & 10.5 & 11.2 & 7.7 & 10.5 & 10.5 & 13.7 & 14.8 & 14.7 & 15.2 \\
$\mathrm{Ray}(1)$ & 100.0 & 100.0 & 100.0 & 100.0 & 100.0 & 100.0 & 100.0 & 100.0 & 100.0 & 100.0 & 100.0 & 100.0 & 100.0 & 100.0 & 100.0 & 100.0 & 100.0 & 100.0 \\
$\mathrm{W}(0.5)$ & 80.4 & 71.6 & 76.0 & 91.8 & 60.5 & 62.6 & 81.6 & 65.0 & 54.6 & 73.2 & 64.5 & 64.9 & 81.8 & 64.4 & 76.0 & 66.4 & 58.1 & 59.2 \\
$\mathrm{W}(1.2)$ & 96.4 & 90.3 & 96.7 & 96.3 & 98.2 & 98.0 & 94.6 & 98.2 & 98.4 & 97.5 & 97.5 & 95.0 & 96.3 & 98.5 & 97.1 & 98.1 & 98.3 & 98.5 \\
$\mathrm{W}(1.75)$ & 100.0 & 100.0 & 100.0 & 100.0 & 100.0 & 100.0 & 100.0 & 100.0 & 100.0 & 100.0 & 100.0 & 100.0 & 100.0 & 100.0 & 100.0 & 100.0 & 100.0 & 100.0 \\
$\mathrm{L}(2)$ & 57.7 & 56.4 & 57.7 & 70.8 & 73.0 & 47.1 & 68.3 & 58.0 & 52.0 & 56.7 & 50.6 & 37.9 & 67.0 & 38.4 & 64.9 & 84.3 & 78.6 & 78.6 \\
$\mathrm{Burr}(1.5,0.5,0.5)$ & 46.8 & 42.0 & 42.2 & 76.1 & 68.6 & 28.6 & 49.4 & 19.2 & 19.5 & 23.2 & 23.5 & 19.4 & 55.4 & 23.8 & 54.6 & 77.1 & 74.0 & 74.4 \\
$\mathrm{IG}(1,1)$ & 99.9 & 99.8 & 100.0 & 100.0 & 99.9 & 99.9 & 100.0 & 99.9 & 100.0 & 99.9 & 99.9 & 98.5 & 100.0 & 99.2 & 100.0 & 99.9 & 99.8 & 99.9 \\
$\mathrm{LogW}(0.1)$ & 7.2 & 6.8 & 7.3 & 6.6 & 7.7 & 7.4 & 7.1 & 7.6 & 7.5 & 6.6 & 6.3 & 6.7 & 7.1 & 7.0 & 7.6 & 7.5 & 7.3 & 7.8 \\
$\mathrm{LogW}(0.3)$ & 22.1 & 18.8 & 21.6 & 21.0 & 24.1 & 22.4 & 22.9 & 22.3 & 22.3 & 19.9 & 19.2 & 14.1 & 20.8 & 17.5 & 23.7 & 24.4 & 23.6 & 24.5 \\
$\mathrm{LogW}(1.75)$ & 99.9 & 99.4 & 99.9 & 100.0 & 100.0 & 99.9 & 99.9 & 99.8 & 99.8 & 99.7 & 99.7 & 97.5 & 99.8 & 99.0 & 99.9 & 100.0 & 99.9 & 99.9 \\
$\mathrm{FR}(1,1)$ & 93.9 & 95.7 & 95.4 & 98.5 & 91.7 & 84.0 & 97.9 & 92.0 & 86.6 & 91.5 & 85.2 & 82.7 & 97.9 & 76.7 & 94.7 & 95.7 & 92.9 & 93.0 \\
$\mathrm{HN}(0.8)$ & 94.2 & 88.6 & 95.6 & 94.9 & 97.6 & 96.9 & 91.3 & 97.5 & 97.4 & 96.2 & 95.9 & 97.2 & 94.4 & 98.1 & 93.8 & 97.6 & 97.6 & 97.9 \\
$\mathrm{HN}(1)$ & 97.0 & 93.0 & 97.9 & 97.4 & 98.7 & 98.5 & 94.9 & 98.9 & 98.9 & 98.1 & 98.0 & 98.3 & 96.9 & 99.3 & 96.5 & 98.6 & 98.9 & 99.0 \\
$\chi^2_4$ & 100.0 & 100.0 & 100.0 & 100.0 & 100.0 & 100.0 & 100.0 & 100.0 & 100.0 & 100.0 & 100.0 & 100.0 & 100.0 & 100.0 & 100.0 & 100.0 & 100.0 & 100.0 \\
$\mathrm{DH}(0.2)$ & 28.6 & 23.7 & 29.1 & 27.4 & 32.9 & 31.6 & 28.7 & 30.3 & 31.4 & 26.3 & 26.9 & 18.0 & 26.1 & 24.3 & 31.8 & 32.9 & 32.5 & 33.3 \\
$\mathrm{DH}(0.4)$ & 74.9 & 63.7 & 74.7 & 74.2 & 80.8 & 79.7 & 74.1 & 74.6 & 77.7 & 70.7 & 72.9 & 51.2 & 70.7 & 66.2 & 79.4 & 80.4 & 80.4 & 80.9 \\
$\mathrm{DH}(0.7)$ & 99.2 & 96.9 & 99.2 & 99.2 & 99.7 & 99.6 & 99.1 & 99.2 & 99.5 & 98.9 & 99.3 & 93.9 & 98.8 & 98.3 & 99.6 & 99.6 & 99.7 & 99.7 \\
$\mathrm{LLogis}(2,1)$ & 100.0 & 99.7 & 99.9 & 99.9 & 99.9 & 99.9 & 100.0 & 100.0 & 99.9 & 99.9 & 99.9 & 99.4 & 100.0 & 99.7 & 100.0 & 99.9 & 99.9 & 99.9 \\
$\mathrm{LFR}(0.2)$ & 76.9 & 68.4 & 81.0 & 78.4 & 85.5 & 84.6 & 72.0 & 86.7 & 86.9 & 82.3 & 82.1 & 79.4 & 77.5 & 88.5 & 78.4 & 84.7 & 85.5 & 86.7 \\
$\mathrm{LFR}(0.5)$ & 86.5 & 79.1 & 89.5 & 87.6 & 93.2 & 92.1 & 82.0 & 93.5 & 93.5 & 90.6 & 90.1 & 90.5 & 87.2 & 94.8 & 86.9 & 92.9 & 93.0 & 93.8 \\
$\mathrm{LFR}(0.8)$ & 90.7 & 83.8 & 92.8 & 91.3 & 95.6 & 94.8 & 86.9 & 95.8 & 95.6 & 93.6 & 93.2 & 94.5 & 91.3 & 96.8 & 91.0 & 95.5 & 95.6 & 96.3 \\
$F_{1,3}$ & 55.3 & 97.9 & 92.3 & 90.8 & 11.7 & 46.5 & 35.9 & 91.2 & 90.1 & 91.5 & 90.6 & 75.2 & 71.1 & 90.9 & 77.5 & 4.8 & 5.2 & 5.2 \\
\hline
\end{tabular}
}
\end{sidewaystable}

\begin{sidewaystable}[ht]
\caption{Empirical size and power (percentages, rounded to one decimal place) of the proposed and competing tests ($n=100$, $N=10{,}000$, $B=500$, $\eta=0.05$). The distribution $F_{1,3}$ is the alternative in Proposition~\ref{prop:continuous_invisible} with $\lambda=1$ and $b=3$.} \label{power100}
\resizebox{\linewidth}{!}{
\begin{tabular}{|l|llllllllllllllllll|} \hline
Distribution  & $O_n$ & $KS_n$ & $CM_n$ & $AD_n$ & $G_{n,2}$ & $M_{n,2}$ & $|I_{n,2}|$ & $S_{n,3,2}^{(1)}$ & $S_{n,3,2}^{(2)}$ & $T_{n,3,2}^{(1)}$ & $T_{n,3,2}^{(2)}$ & $N^{(1)}_{n,2,0.5}$ & $N^{(1)}_{n,3,0.5}$ & $N^{(2)}_{n,2,0.5}$ & $N^{(2)}_{n,3,0.5}$ & ${T}_n$ & $S_{n,K}$ & $S_{n,B}$ \\ \hline
$\mathcal P(0.5)$ & 5.4 & 4.7 & 5.1 & 5.1 & 5.5 & 5.3 & 5.3 & 4.9 & 4.8 & 4.8 & 4.8 & 5.2 & 5.0 & 5.0 & 5.0 & 5.5 & 5.2 & 5.3 \\
$\mathcal P(1)$ & 5.2 & 4.8 & 4.9 & 5.0 & 5.3 & 5.4 & 5.4 & 4.7 & 4.5 & 4.8 & 4.6 & 4.9 & 4.8 & 4.9 & 4.9 & 5.2 & 5.1 & 5.3 \\
$\mathcal P(2)$ & 4.9 & 5.0 & 4.7 & 4.8 & 5.0 & 4.9 & 5.0 & 4.3 & 4.3 & 4.4 & 4.5 & 4.8 & 4.7 & 4.8 & 4.9 & 5.0 & 4.8 & 5.0 \\
\hline
$\Gamma(0.5)$ & 86.7 & 80.5 & 84.6 & 96.5 & 41.6 & 67.8 & 88.9 & 68.5 & 59.3 & 75.4 & 68.0 & 74.6 & 84.6 & 55.6 & 80.1 & 42.1 & 35.2 & 36.4 \\
$\Gamma(1.2)$ & 99.9 & 99.4 & 100.0 & 100.0 & 100.0 & 100.0 & 99.8 & 100.0 & 100.0 & 100.0 & 100.0 & 99.8 & 99.9 & 100.0 & 100.0 & 100.0 & 100.0 & 100.0 \\
$\Gamma(2)$ & 100.0 & 100.0 & 100.0 & 100.0 & 100.0 & 100.0 & 100.0 & 100.0 & 100.0 & 100.0 & 100.0 & 100.0 & 100.0 & 100.0 & 100.0 & 100.0 & 100.0 & 100.0 \\
$\mathrm{InvB}(0.1)$ & 9.0 & 8.0 & 8.9 & 8.5 & 9.4 & 8.9 & 9.4 & 9.7 & 9.5 & 8.5 & 8.2 & 7.4 & 8.7 & 8.0 & 9.4 & 9.5 & 9.2 & 9.5 \\
$\mathrm{InvB}(2)$ & 100.0 & 100.0 & 100.0 & 100.0 & 100.0 & 100.0 & 100.0 & 100.0 & 100.0 & 100.0 & 100.0 & 99.9 & 100.0 & 100.0 & 100.0 & 100.0 & 100.0 & 100.0 \\
$\mathrm{Tilt}(1)$ & 45.0 & 39.4 & 47.4 & 44.4 & 49.3 & 50.8 & 43.5 & 40.6 & 43.9 & 37.5 & 40.8 & 20.4 & 32.8 & 32.4 & 44.2 & 45.0 & 48.1 & 48.1 \\
$\mathrm{Tilt}(2)$ & 85.0 & 78.5 & 86.7 & 84.9 & 87.1 & 89.4 & 83.7 & 75.4 & 79.3 & 73.0 & 77.2 & 43.4 & 66.5 & 63.6 & 81.4 & 82.0 & 85.7 & 85.6 \\
$\mathrm{Tilt}(3)$ & 97.0 & 94.2 & 97.7 & 96.9 & 97.2 & 98.4 & 96.4 & 91.6 & 93.7 & 90.4 & 92.6 & 61.8 & 85.6 & 82.4 & 94.5 & 94.2 & 96.3 & 96.3 \\
$\mathrm{Ben}(0.1)$ & 12.9 & 12.7 & 14.9 & 13.4 & 16.8 & 17.3 & 12.6 & 15.7 & 17.6 & 13.6 & 14.9 & 9.8 & 11.8 & 14.1 & 15.2 & 15.7 & 16.5 & 16.9 \\
$\mathrm{Ben}(0.3)$ & 42.1 & 39.3 & 47.4 & 44.4 & 54.2 & 54.1 & 40.0 & 51.6 & 56.2 & 46.9 & 50.5 & 33.6 & 41.5 & 49.6 & 49.1 & 51.4 & 54.2 & 54.6 \\
$\mathrm{Ben}(0.7)$ & 79.4 & 73.3 & 83.7 & 81.9 & 89.0 & 88.4 & 77.0 & 88.4 & 90.2 & 85.8 & 87.3 & 76.2 & 81.2 & 88.5 & 84.9 & 88.1 & 89.5 & 89.5 \\
$\mathrm{LogG}(0.5)$ & 78.6 & 67.2 & 75.9 & 79.6 & 83.1 & 78.0 & 79.9 & 65.3 & 66.7 & 62.9 & 64.3 & 38.7 & 66.2 & 52.4 & 76.1 & 84.1 & 83.2 & 83.4 \\
$\mathrm{LogG}(1.2)$ & 100.0 & 99.9 & 100.0 & 100.0 & 100.0 & 100.0 & 100.0 & 99.6 & 99.5 & 99.4 & 99.4 & 91.5 & 99.6 & 97.6 & 99.9 & 100.0 & 100.0 & 100.0 \\
$\mathrm{LogG}(2)$ & 100.0 & 100.0 & 100.0 & 100.0 & 100.0 & 100.0 & 100.0 & 100.0 & 100.0 & 100.0 & 100.0 & 99.4 & 100.0 & 100.0 & 100.0 & 100.0 & 100.0 & 100.0 \\
$\mathrm{LN}(1)$ & 100.0 & 100.0 & 100.0 & 100.0 & 100.0 & 100.0 & 100.0 & 100.0 & 100.0 & 100.0 & 100.0 & 99.9 & 100.0 & 100.0 & 100.0 & 100.0 & 100.0 & 100.0 \\
$\mathrm{LN}(1.2)$ & 96.5 & 90.6 & 96.0 & 97.9 & 97.9 & 97.3 & 97.0 & 94.4 & 95.8 & 93.2 & 94.7 & 78.9 & 95.5 & 89.6 & 97.5 & 98.0 & 97.6 & 97.6 \\
$\mathrm{LN}(1.5)$ & 21.4 & 16.5 & 20.9 & 23.9 & 29.3 & 27.5 & 21.5 & 19.3 & 23.0 & 16.5 & 19.7 & 10.4 & 17.4 & 16.2 & 24.6 & 27.7 & 27.7 & 28.2 \\
$\mathrm{Ray}(1)$ & 100.0 & 100.0 & 100.0 & 100.0 & 100.0 & 100.0 & 100.0 & 100.0 & 100.0 & 100.0 & 100.0 & 100.0 & 100.0 & 100.0 & 100.0 & 100.0 & 100.0 & 100.0 \\
$\mathrm{W}(0.5)$ & 97.1 & 95.3 & 96.7 & 99.6 & 87.8 & 90.1 & 97.6 & 95.1 & 91.2 & 96.3 & 93.4 & 92.9 & 98.0 & 91.6 & 96.1 & 91.1 & 86.7 & 87.0 \\
$\mathrm{W}(1.2)$ & 100.0 & 99.8 & 100.0 & 100.0 & 100.0 & 100.0 & 99.9 & 100.0 & 100.0 & 100.0 & 100.0 & 100.0 & 100.0 & 100.0 & 100.0 & 100.0 & 100.0 & 100.0 \\
$\mathrm{W}(1.75)$ & 100.0 & 100.0 & 100.0 & 100.0 & 100.0 & 100.0 & 100.0 & 100.0 & 100.0 & 100.0 & 100.0 & 100.0 & 100.0 & 100.0 & 100.0 & 100.0 & 100.0 & 100.0 \\
$\mathrm{L}(2)$ & 91.7 & 94.6 & 93.1 & 99.0 & 97.3 & 78.5 & 96.8 & 89.0 & 82.7 & 88.3 & 81.9 & 74.0 & 97.3 & 72.1 & 95.3 & 99.3 & 98.6 & 98.6 \\
$\mathrm{Burr}(1.5,0.5,0.5)$ & 70.0 & 74.3 & 73.3 & 95.7 & 91.8 & 49.5 & 73.9 & 54.2 & 51.5 & 57.7 & 54.8 & 49.2 & 85.5 & 52.8 & 84.0 & 95.6 & 94.4 & 94.4 \\
$\mathrm{IG}(1,1)$ & 100.0 & 100.0 & 100.0 & 100.0 & 100.0 & 100.0 & 100.0 & 100.0 & 100.0 & 100.0 & 100.0 & 100.0 & 100.0 & 100.0 & 100.0 & 100.0 & 100.0 & 100.0 \\
$\mathrm{LogW}(0.1)$ & 9.2 & 8.6 & 9.1 & 8.9 & 9.8 & 9.4 & 10.1 & 9.8 & 9.6 & 8.8 & 8.1 & 7.4 & 9.3 & 7.9 & 10.1 & 10.1 & 9.7 & 10.0 \\
$\mathrm{LogW}(0.3)$ & 42.0 & 33.0 & 39.2 & 41.6 & 43.6 & 40.8 & 43.6 & 38.5 & 39.2 & 35.7 & 35.7 & 23.0 & 37.9 & 30.6 & 43.5 & 44.8 & 43.5 & 43.8 \\
$\mathrm{LogW}(1.75)$ & 100.0 & 100.0 & 100.0 & 100.0 & 100.0 & 100.0 & 100.0 & 100.0 & 100.0 & 100.0 & 100.0 & 100.0 & 100.0 & 100.0 & 100.0 & 100.0 & 100.0 & 100.0 \\
$\mathrm{FR}(1,1)$ & 100.0 & 100.0 & 100.0 & 100.0 & 99.9 & 99.1 & 100.0 & 99.9 & 99.7 & 99.9 & 99.6 & 99.7 & 100.0 & 98.7 & 100.0 & 100.0 & 99.9 & 99.9 \\
$\mathrm{HN}(0.8)$ & 100.0 & 99.7 & 100.0 & 100.0 & 100.0 & 100.0 & 99.9 & 100.0 & 100.0 & 100.0 & 100.0 & 100.0 & 100.0 & 100.0 & 99.9 & 100.0 & 100.0 & 100.0 \\
$\mathrm{HN}(1)$ & 100.0 & 99.9 & 100.0 & 100.0 & 100.0 & 100.0 & 99.9 & 100.0 & 100.0 & 100.0 & 100.0 & 100.0 & 100.0 & 100.0 & 100.0 & 100.0 & 100.0 & 100.0 \\
$\chi^2_4$ & 100.0 & 100.0 & 100.0 & 100.0 & 100.0 & 100.0 & 100.0 & 100.0 & 100.0 & 100.0 & 100.0 & 100.0 & 100.0 & 100.0 & 100.0 & 100.0 & 100.0 & 100.0 \\
$\mathrm{DH}(0.2)$ & 53.0 & 43.7 & 52.0 & 52.7 & 59.1 & 57.3 & 53.6 & 50.5 & 53.9 & 47.2 & 49.6 & 30.9 & 47.7 & 43.3 & 56.7 & 58.7 & 58.0 & 58.9 \\
$\mathrm{DH}(0.4)$ & 97.2 & 92.3 & 97.1 & 97.6 & 98.6 & 98.4 & 97.2 & 96.4 & 97.5 & 95.4 & 96.6 & 84.9 & 95.8 & 94.0 & 98.1 & 98.4 & 98.5 & 98.5 \\
$\mathrm{DH}(0.7)$ & 100.0 & 100.0 & 100.0 & 100.0 & 100.0 & 100.0 & 100.0 & 100.0 & 100.0 & 100.0 & 100.0 & 100.0 & 100.0 & 100.0 & 100.0 & 100.0 & 100.0 & 100.0 \\
$\mathrm{LLogis}(2,1)$ & 100.0 & 100.0 & 100.0 & 100.0 & 100.0 & 100.0 & 100.0 & 100.0 & 100.0 & 100.0 & 100.0 & 100.0 & 100.0 & 100.0 & 100.0 & 100.0 & 100.0 & 100.0 \\
$\mathrm{LFR}(0.2)$ & 97.2 & 94.9 & 98.6 & 98.5 & 99.4 & 99.3 & 95.4 & 99.4 & 99.4 & 99.0 & 99.1 & 98.7 & 97.8 & 99.7 & 97.7 & 99.3 & 99.5 & 99.5 \\
$\mathrm{LFR}(0.5)$ & 99.3 & 98.2 & 99.7 & 99.7 & 99.9 & 99.8 & 98.7 & 99.9 & 99.9 & 99.8 & 99.8 & 99.8 & 99.6 & 99.9 & 99.4 & 99.8 & 99.9 & 99.9 \\
$\mathrm{LFR}(0.8)$ & 99.8 & 99.1 & 99.9 & 99.9 & 100.0 & 99.9 & 99.4 & 100.0 & 100.0 & 99.9 & 99.9 & 99.9 & 99.9 & 100.0 & 99.8 & 100.0 & 100.0 & 100.0 \\
$F_{1,3}$ & 83.3 & 100.0 & 99.9 & 99.9 & 20.8 & 80.5 & 58.4 & 99.8 & 99.8 & 99.8 & 99.8 & 98.8 & 97.6 & 99.9 & 98.3 & 4.5 & 5.8 & 6.0 \\
\hline
\end{tabular}
}
\end{sidewaystable}

\end{appendix}

\end{document}